\documentclass[11pt]{amsart}
\usepackage[latin1]{inputenc}
\usepackage{amsmath}
\usepackage{amssymb}
\usepackage[all]{xy}
\usepackage{mathrsfs}
\usepackage{amsthm}
\usepackage{graphicx}

\theoremstyle{plain}
\newtheorem{thm}[equation]{Theorem}
\newtheorem{lem}[equation]{Lemma}
\newtheorem{prop}[equation]{Proposition}
\newtheorem{cor}[equation]{Corollary}

\theoremstyle{definition}
\newtheorem{defin}[equation]{Definition}
\newtheorem{remark}[equation]{Remark}

\numberwithin{equation}{subsection}

\def\sheafEnd{\mathcal{E} \hspace{-1pt} \mathit{nd}}
\def\sheafHom{\mathcal{H} \hspace{-1pt} \mathit{om}}

\newcommand{\nc}{\newcommand}
\newcommand{\rc}{\renewcommand}

\nc{\Aut}{{	\operatorname{Aut}	}}
\nc{\codim}{{	\operatorname{codim}	}}
\nc{\Ob}{{	\operatorname{Ob}	}}
\nc{\PGL}{{	\operatorname{PGL}	}}
\nc{\supp}{{	\operatorname{supp}	}}
\nc{\tr}{{	\operatorname{tr}	}}

\nc{\Rep}{{	{\cal{R}}ep		}}
\nc{\one}{{	\mbox{\bf{1}}		}}
\nc{\iso}{	\overset{\sim}{\lra}	}

\nc{\nen}{\newenvironment}

\nc{\pr}{\protect}

\nc{\nn}{{\newline}}
\nc{\np}{{\newpage}}	

\nc{\lab}{	\label}
\nc{\npp}{{	\newpage\setcounter{page}{0}	}}
\nc{\setpa}{		\setcounter{part}		}
\nc{\setse}{		\setcounter{section}	}
\nc{\setsus}{		\setcounter{subsection}		}
\nc{\setsss}{		\setcounter{subsubsection}	}
\nc{\setpage}{		\setcounter{page}	}

\nc{\nfd}{ $$\text{ This version is preliminary and approximate, 
		             it is not for distribution. }$$	}

\nc{\noi}{{\noindent}}
\nc{\pf}{{	\noindent {\em Proof.}		}}
					
\nc{\epf}{ \fbox{\bf QED}	}

\nc{\heart}{{\tiny \cen{\tiny $\heartsuit $ }	}} 
\nc{\cont}{\tableofcontents}

\nc{\sbr}{{	\smallpagebreak	}}
\nc{\mbr}{{	\medpagebreak	}}
\nc{\bbr}{{	\bigpagebreak	}}
\nc{\bbb}{ 	\boldsymbol 	}  

\nc{\bib}{		}
\nc{\bit}[1]{	\bibitem[#1]{#1} 		}

\rc{\b}{ 	\big         			}  
\nc{\lam}[1]{{ 	\text{\large $#1$	}	}}  
\nc{\smm}[1]{{ 	\text{\small $#1$	}	}}  
\nc{\fom}[1]{{ 	\text{\footnotesize $#1$	}	}}  
\nc{\tinm}[1]{{ \text{\tiny $#1$	}	}}  

\nc{\bu}{ \bullet         }  			
\nc{\bbu}{ \aa{\bbb \bullet}         }  	
\nc{\bus}{{	^\bullet	}}	 	
\nc{\bui}{{	_\bullet	}}	 	

\nc{\bem}{{	\begin{em}	}}
\nc{\eem}{{	\end{em} 	}}

\nc{\bbox}{{	\blackbox	}}	
\nc{\bx}{	\boxed	}		
\nc{\tbx}[1]{{\boxed{\tx{#1}}}}		

\nc{\mmbox}[1]{{	\mbox{$#1$}	}}	
\nc{\tbox}[1]{{		\mbox{\tx{#1}}	}}
 	
\nc{\ot}{		\leftarrow			}
\nc{\tto}{		\longrightarrow			}
\nc{\ott}{		\longleftarrow			}
\nc{\too}[1]{{		\aa{#1}\rightarrow			}}
\nc{\oot}[1]{{		\aa{#1}\leftarrow			}}
\nc{\ttoo}[1]{{		\aa{#1}\longrightarrow			}}
\nc{\oott}[1]{{		\aa{#1}\longleftarrow			}}

\nc{\Too}[2]{{		\aa{#1}{\bb{#2}\rightarrow}		}}
\nc{\ooT}[2]{{		\aa{#1}{\bbb{#2}\leftarrow}		}}
\nc{\TToo}[2]{{		\aa{#1}{\bb{#2}\longrightarrow}		}}
\nc{\ooTT}[2]{{		\aa{#1}{\bbb{#2}\longleftarrow}		}}
\nc{\toot}[2]{{		\aa{#1}{\bb{#2}\rightleftarrows}	}}
\nc{\ttoot}[2]{{	\aa{#1}{\bb{#2}\rightleftrightarrows}	}}

\nc{\ra}{{	\rightarrow		}}
\nc{\laa}{{	\leftarrow	}}	
\nc{\lra}{{\longrightarrow}}

\nc{\lr}{{\leftrightarrow}}     	
\nc{\lrs}{{\rightleftarrows}}     	

\nc{\imp}{{\Rightarrow}}        	
\nc{\impp}{{\Leftarrow}}        	
\nc{\eq}{{\Leftrightarrow}}        	
\nc{\impl}{{\Longrightarrow}}        	
\nc{\imppl}{{\Longleftarrow}}        	
\nc{\eql}{{\Longleftrightarrow}}        	
	\nc{\Ra}{{\Rightarrow}}         	
	\nc{\LRa}{{\Leftrightarrow}}        	

\nc{\inj}{{\pr	\hookrightarrow	}}    		
\nc{\injj}{{\pr	\hookleftarrow	}}    		

\nc{\sur}{{	\twoheadrightarrow	}}	
\nc{\surr}{{	\twoheadleftarrow	}}	
			
\nc{\mm}{{	\mapsto		}}     		
\nc{\mmm}{{	\leftarrow\shortmid }}		
\nc{\ainj}[1]{{\aa{#1}{\pr\hookrightarrow}	}}    	
\nc{\ainjj}[1]{{\aa{#1}{\pr\hookleftarrow}	}}    	

\nc{\asur}[1]{{	\aa{#1}\twoheadrightarrow	}}	
\nc{\asurr}[1]{{\aa{#1}\twoheadleftarrow	}}	
			
\nc{\amm}[1]{{	\aa{#1}\mapsto		}}     	
\nc{\ammm}[1]{{	\aa{#1}\leftarrow\shortmid }}	

\nc{\va}{{\uparrow}}              		

\nc{\syp}[1]{	^{ (#1) }		} 	
\nc{\up}[1]{	^{ (#1) }		} 	
\nc{\lp}[1]{	_{ (#1) }		}	
\nc{\hp}[1]{	^{ [#1] }		}	
\nc{\cle}{\preceq}		
\nc{\cl}{\prec}			
\nc{\cge}{\succeq}		
\nc{\cg}{\succ}			

\nc{\bb}{	\pr\underset 	}           
\rc{\aa}{ 	\pr\overset 	}            

\nc{\indd}{{ ${} \ \ \ \ \  \ \        {} $	}}	
\nc{\inddd}{{ 	\indd\indd			}}	
\nc{\nnd}{{ 	\nn  \indd 			}}	
\nc{\nndb}{{ 	\nn  \indd $\bullet$		}}	
		
\nc{\bce}{	\begin{center}	}
\nc{\ece}{	\end  {center}	}
\nc{\cen}[1]{	\begin{center}	{\em  #1}	\end  {center}	}

\nc{\bss}{{\backslash}}           		
\nc{\barr}{ 	\overline 	}      		
\nc{\ud}{	\underline	}		

\nc{\ti}{\tilde}              
\nc{\tii}{\widetilde}         

\nc{\hatt}{\widehat}				
\nc{\hata}{{	\bbb{ \hat{} }		}}	

\nc{\ch}{\check}              			
\nc{\cha}{{ 	\bbb{ \check{} }	}}      

\nc{\sub}{{	\subseteq	}}         
\nc{\subb}{{	\supseteq	}}         
\nc{\nsub}{{	\nsubseteq	}}         
\nc{\nsubb}{{	\nsupseteq	}}         %

\nc{\nin}{{	\notin	}}

\nc{\lb}{\langle}             				
\nc{\rb}{\rangle}

\nc{\lB}{	\left(	}             			
\nc{\rB}{	\right)	}
\nc{\BBl}{{	\bbb{ \left( \right.}	}}             	
\nc{\BBr}{{	\bbb{ \left. \right)}	}}

\nc{\Pa}[2]{ {\lb} #1 {,} #2 {\rb} }				

\nc{\cD}[1]{ \tx{ $$\CD {#1} \endCD $$ }  }		

\nc{\mat} {		\left(		\matrix	}	
\nc{\emat}{		\endmatrix	\right)	}
\nc{\sm} {		\left(		\smallmatrix	}	
\nc{\esm}{		\endsmallmatrix	\right)	}
\nc{\smat} {		\left(		\smallmatrix	}	
\nc{\esmat}{		\endsmallmatrix	\right)	}

\nc{\matr} {		\left[		\matrix	}	
\nc{\ematr}{		\endmatrix	\right]	}
\nc{\smr} {		\left[		\smallmatrix	}	
\nc{\esmr}{		\endsmallmatrix	\right]	}
\nc{\smatr} {		\left[		\smallmatrix	}	
\nc{\esmatr}{		\endsmallmatrix	\right]	}

\nc{\imat} {		\left.		\matrix	}	
\nc{\eimat}{		\endmatrix	\right.	}
\nc{\ism} {		\left.		\smallmatrix	}	
\nc{\eism}{		\endsmallmatrix	\right.	}

\nc{\ca}{		\left\{		\smallmatrix	}	
\nc{\eca}{		\endsmallmatrix	\right\}	}
\nc{\Ca}{		\left\{		\matrix		}	
\nc{\Eca}{		\endmatrix	\right.		}	
\nc{\eCa}{		\endmatrix	\right\}	}	

\nc{\com}{	\begin{diagram}	}
\nc{\ecom}{	  \end{diagram}	}

\nc{\tab}{	\begin{tabular}		}
\nc{\etab}{	\end{tabular}		}	
\nc{\hl}{{	\hline			}}

\nc{\Eq}{	\begin{equation}	}
\nc{\Eeq}{	\end{equation}	}
\nc{\aln}{	\begin{align}	}
\nc{\ealn}{	\end{align}	}

\nc{\Rpart}{	\rc{\thepart}{\Roman{part}}	}
\nc{\Apart}{	\rc{\thepart}{\arabic{part}}	}
		
\nc{\rref}[2]{\ref{#1}.\ref{#2}}

\nc{\pa}[1]{ 	\part{#1}		}

\nc{\se}[1]{ 	\section{\bf#1}		}
\nc{\ses}[1]{ 	\section*{\bf#1}		}
\nc{\sus}{ 	\subsection		}
\nc{\sss}{ 	\subsubsection		}

\nc{\Lem}{ 	\subsection{Lemma}		}
\nc{\slem}{ 	\subsubsection*{Lemma}		}
\nc{\sublem}{ 	\subsubsection{ Sublemma}	}
\nc{\ssublem}{ \subsubsection*{ Sublemma}	}

\nc{\Lemm}{ 	\subsection{Lemma}		}
\nc{\lemm}{ 	\subsubsection{Lemma}		}
\nc{\slemm}{ 	\subsubsection*{Lemma}		}
\nc{\sublemm}{ 	\subsubsection{ Sublemma}	}
\nc{\ssublemm}{ \subsubsection*{ Sublemma}	}

\nc{\Pro}{ 	\subsection{Proposition}	}
\nc{\pro}{ 	\subsubsection{Proposition}	}
\nc{\spro}{ 	\subsubsection*{Proposition}	}

\nc{\Cor}{ 	\subsection{Corollary}		}
\nc{\scor}{ 	\subsubsection*{Corollary}	}

\nc{\Corr}{ 	\subsection{Corollary}		}
\nc{\corr}{ 	\subsubsection{Corollary}	}
\nc{\scorr}{ 	\subsubsection*{Corollary}	}

\nc{\Theo}{ 	\subsection{Theorem}		}		
\nc{\theo}{ 	\subsubsection{Theorem}		}
\nc{\stheo}{ 	\subsubsection*{Theorem}	}

\nc{\pretheo}{ 	\subsubsection{Pretheorem}	}

\nc{\rem}{ 	\subsubsection{Remark}		}
\nc{\srem}{ 	\subsubsection*{Remark}	}

\nc{\rems}{ 	\subsubsection{Remarks}		}

\nc{\srems}{ 	\subsubsection*{Remarks}	}

\nc{\Def}{ 	\subsection{Definition}		}
\nc{\ddef}{ 	\subsubsection{Definition}	}

\nc{\comm}{ 	\subsubsection{Comment}		}
\nc{\scomm}{ 	\subsubsection*{Comment}	}
\nc{\comms}{ 	\subsubsection{Comments}		}
\nc{\scomms}{ 	\subsubsection*{Comments}	}

\nc{\claim}{ 	\subsubsection{Claim}		}
\nc{\sclaim}{ 	\subsubsection*{Claim}	}

\nc{\nota}{ 	\subsubsection{Notation}	}

\nc{\conj}{ 	\subsubsection{Conjecture}	}
\nc{\sconj}{ 	\subsubsection*{Conjecture}	}

\nc{\ex}{ 	\subsubsection{Example}		}
\nc{\sex}{ 	\subsubsection*{Example}	}
\nc{\exs}{ 	\subsubsection{Examples}	}
\nc{\sexs}{ 	\subsubsection*{Examples}	}
\nc{\Ex}{ 	\subsection{Example}		}
\nc{\sEx}{ 	\subsection*{Example}	}
\nc{\Exs}{ 	\subsection{Examples}	}
\nc{\sExs}{ 	\subsection*{Examples}	}

\nc{\que}{ 	\subsubsection{Question}	}
\nc{\ques}{ 	\subsubsection{Questions}	}
\nc{\sque}{ 	\subsubsection*{Question}	}
\nc{\sques}{ 	\subsubsection*{Questions}	}

\nc{\bi}{	\begin{itemize}\item		}
\rc{\i}{	\item			}
\nc{\ei}{ \end{itemize}	} 
\nc{\ben}{	\begin{enumerate}\item		}
\nc{\een}{	\end{enumerate}			}
\nc{\ftt}[1]{{\footnote{#1}}}
\nc{\fttt}[1]{{$^($\footnote{#1}$^)$}}
\nc{\bftt}[1]{\footnote{#1}}
\nc{\f}[1]{ \fbox{$ $}\footnote{ \fbox{!}#1 }\fbox{$ $}		}

\nc{\Ao}{{	\A^1	}}
\nc{\Po}{{	\P^1	}}

\nc{\h}{{	\hslash	}}	

\nc{\All}{{	\forall		}}
\nc{\Exx}{{	\exists 	}}
\nc{\yy}{\infty}                       
\nc{\ys}{{  \frac{\infty}{2}  }}

\nc{\ii}{{i\in I}}
\nc{\ww}{{w\in W}}
\nc{\SES}[5]{{	0 @>>> {#1} @>{#2}>> {#3} @>{#4}>> {#5} @>>> 0	}}
\nc{\Ses}[3] {{	0 @>>> {#1} @>>>     {#2} @>>>     {#3} @>>> 0	}}
\nc{\pl}{{\oplus}}              		
\nc{\tim}{{\times}}             
\nc{\btim}{{\boxtimes}}
\nc{\ltim}{\ltimes}                  	%
\nc{\rtim}{\rtimes}			%
\nc{\ltr}{\triangleleft}        %
\nc{\rtr}{\triangleright}       %

\nc{\ten}{{	\otimes		}}            
\nc{\Lten}{{	\aa{L}\otimes	}}            
\nc{\Ltim}{{	\aa{L}\times	}}            
\nc{\Rcap}{{	\aa{R}\cap	}}            
\nc{\tenA}{	\bb{A}\ten	}
\nc{\tenB}{	\bb{B}\ten	}
\nc{\tenZ}{	\bb{\Z}\ten	}
\nc{\tenR}{	\bb{\R}\ten	}
\nc{\tenC}{	\bb{\C}\ten	}
\nc{\tenk}{	\bb{\k}\ten	}
\nc{\bten}{{\boxtimes}}         		

\nc{\con}{{ @>{\protect\cong}>> }}  	
\nc{\conl}{{ 	@>{\cong}>>	}}  	
\nc{\conn}{{    @<{\cong}<<  	}}  	
\nc{\Con}{{	\equiv		}}	
\nc{\appr}{{	\sim		}}	
\nc{\eqr}{{	\sim		}}	
\nc{\equi}{{	\sim		}}	

\nc{\fra}{ 	\frac	}     	
\nc{\ffr}[2]{{ 	\text{\footnotesize $\frac{#1}{#2}$	}	}}  

\nc{\ha}{{ \frac{1}{2} }}     		
	\nc{\half}{{ \frac{1}{2} }}    	

\nc{\ci}{{\circ}}               
\nc{\cd }{{\cdot}}            	
\nc{\cddd}{{\cdot\cdot\cdot}}	

\nc{\ox}{{	\OO_X		}}               
\nc{\omx}{{	\om_X		}}               
\nc{\Omx}{{	\Om_X^1		}}               
\nc{\qcoh}{{	q\CC oh		}}               %

\nc{\xt}{{	X_*(T)		}}
\nc{\Xt}{{	X^*(T)		}}

\nc{\cfm}{{	co\fm		}}	

\nc{\cupp}{\bigcup}             
\nc{\capp}{\bigcap}
\nc{\pll}{\bigoplus}

\nc{\pii}{\prod}                
\nc{\ppii}{\bigprod}            
\nc{\cci}{\sqcup}              
\nc{\ccii}{\bigsqcup}

\nc{\wwe}{\bigwedge}            
\nc{\cce}{\bigcoprod}           

\nc{\aaa}{	\stackerel	}	

\nc{\edd}{{ \end{document}	}} 

\nc{\tx}{	\text		}		
\nc{\df}{{ \protect\overset{ \text{def}}= 	}}		
\nc{\dff}{{ \ \df\				}}		
\nc{\inv}{{ {}^{-1}      }}			

\nc{\thh}{	^{\text{th}}	}                     	
\nc{\st}{	^{\text{st}}	}                     	
\nc{\nd}{	^{\text{nd}}	}                     	
\nc{\rd}{	^{\text{rd}}	}                     	

\nc{\pmo}{{ 	\pm 1		}}
\nc{\mpo}{{ 	\mp 1		}}

\nc{\htt}{  \text{ht}}				

\nc{\emp}{{   \emptyset}}      			

\nc{\cowe}{{	\vee	}}			
\nc{\we}{{\wedge}}				
\nc{\wee}{{	\aa{\bullet}\wedge	}}		
\nc{\wetwo}{{     \pr\overset{2}\wedge       }}	

\nc{\limp}{{	\pr\underset {\leftarrow} \lim		}}	
\nc{\Limp}{{	\pr\underset {\leftarrow} {\bbb\lim}	}}	
\nc{\limi}{{	\pr\underset {\rightarrow}\lim		}}      
\nc{\Limi}{{	\pr\underset {\rightarrow}{\bbb\lim}	}}	
\nc{\llim}[1]{	 \bb{#1}\lim        	}   
\nc{\llimp}[1]{ \bb{#1}{ \pr\underset {\leftarrow} \lim       } }
\nc{\LLimp}[1]{ \bb{#1}{ \pr\underset {\leftarrow} {\bbb\lim} } } 	
\nc{\llimi}[1]{ \bb{#1}{ \pr\underset {\rightarrow}\lim       } }
\nc{\LLimi}[1]{ \bb{#1}{ \pr\underset {\rightarrow}{\bbb\lim} } }	
						
\nc{\ppp}{{ {\Bbb P}^1 }}            		
\nc{\ppn}{{ {\Bbb P}^n }}            		
\nc{\pt}{	{ \text{pt} }	}		
				
\nc{\qlb}{{ \barr{{\Bbb Q}_l} }}      		
\nc{\ffq}{{  {\Bbb F}_q  }}           		
\nc{\ffp}{{  {\Bbb F}_p  }}           		
\nc{\tw}{   {}^{(1)}	}		

\nc{\Ab}{{ 	\AA b 		}}      		%
\nc{\Set}{{ 	\SS et 		}}      		%
\nc{\Top}{{ 	\TT op 		}}      		%
\nc{\Pic}{{ 	\tx{Pic}	}}      		%

\nc{\del}{{\partial }}
\nc{\delb}{{\partial }}
\nc{\dd}[2]{	\fra{d{#1}}{d{#2}}		}
\nc{\ddel}[2]{	\fra{\del{#1}}{\del {#2}}	}
                                        
\nc{\Spec}{{ 	\text{Spec}      		}} 
\nc{\Specf}{{ 	\text{Specf}      		}} 
\nc{\Spf}{{ 	\text{Spf}      		}} 
\nc{\hk}{{     \text{hyperk\"ahler} 	}}
\nc{\susy}{{\text{supersymmetry}}}

\nc{\ie}{{,\ \     \text{i.e.,}\ \ 	}}
\nc{\iif}{{\ \     \text{if}\ \ 	}}
\nc{\aand}{{\ \ \  \text{and}\ \ \ 	}}
\nc{\hence}{{\ \ \ \text{hence}\ \ \ 	}}
\nc{\while}{{\ \ \ \text{while}\ \ \ 	}}
\nc{\with}{{\ \ \  \text{with}\ \ \ 	}}
\nc{\oor}{{\ \     \text{or}\ \ 	}}
\nc{\foor}{{\ \     \text{for}\ \ 	}}
\nc{\suchthat}{{\ \     \text{such that}\ \ 	}}

\nc{\Coker}{{\operatorname{Coker}}}
\rc{\Im}{{ 	\text{Im} 	}}
\nc{\rank}{{	\ \text{rank} 	}}
\nc{\Res}{{	\  \text{Res}   }}
\nc{\End}{{	\text{End}	}}
\nc{\HHom}{{	\text{$\HH$om}	}}
\nc{\RHHom}{{	\text{R$\HH$om} }}
\nc{\RGa}{{	\text{R$\Ga$}	}}
\nc{\EEnd}{{	\text{$\EE nd$}	}}
\nc{\AAut}{{	\text{$\AA ut$}	}}
\nc{\Der}{{	\text{Der}	}}

\nc{\ord	}{{ \text{ord} }}			
\nc{\divv	}{{ \text{div} }}			
\nc{\Lie	}{{ \text{Lie} }}

\nc{\timA} {{   \pr\underset{A}\tim             }}
\nc{\timB} {{   \pr\underset{B}\tim             }}
\nc{\timC} {{   \pr\underset{C}\tim             }}
\nc{\timG} {{   \pr\underset{G}\tim             }}
\nc{\timH} {{   \pr\underset{H}\tim             }}
\nc{\timN} {{   \pr\underset{N}\tim             }}
\nc{\timP}{{    \pr\underset{P}\tim             }}
\nc{\timQ}{{    \pr\underset{Q}\tim             }}
\nc{\timS} {{   \pr\underset{S}\tim             }}
\nc{\timT} {{   \pr\underset{T}\tim             }}
\nc{\timU} {{   \pr\underset{U}\tim             }}
\nc{\timV} {{   \pr\underset{V}\tim             }}
\nc{\timX} {{   \pr\underset{X}\tim             }}
\nc{\timY} {{   \pr\underset{Y}\tim             }}
\nc{\timZ} {{   \pr\underset{Z}\tim             }}

\nc{\ab}{{       ^{\text{ab}}   		}}
\nc{\af}{{       ^{\text{aff}}  		}}

\nc{\cod}{\text{codim}}	

\rc{\AA}{{\mathcal A}}
\nc{\BB}{{\mathcal B}} 
\nc{\CC}{{\mathcal C}}
\nc{\DD}{{\mathcal D}}
\nc{\EE}{{\mathcal E}}
\nc{\FF}{{\mathcal F}}
\nc{\GG}{{\mathcal G}}
\nc{\HH}{{\mathcal H}}
\nc{\II}{{\mathcal I}}
\nc{\JJ}{{\mathcal J}}
\nc{\KK}{{\mathcal K}}
\nc{\LL}{{\mathcal L}}
\nc{\MM}{{\mathcal M}}
\nc{\NN}{{\mathcal N}}
\nc{\OO}{{\mathcal O}}
\nc{\PP}{{\mathcal P}}
\nc{\QQ}{{\mathcal Q}}
\nc{\RR}{{\mathcal R}}
\rc{\SS}{{\mathcal S}}
\nc{\TT}{{\mathcal T}}
\nc{\UU}{{\mathcal U}}
\nc{\VV}{{\mathcal V}}
\nc{\WW}{{\mathcal W}}
\nc{\ZZ}{{\mathcal Z}}
\nc{\XX}{{\mathcal X}}
\nc{\YY}{{\mathcal Y}}

\nc{\A}{{\Bbb A }}
\nc{\B}{{\Bbb B}}
\nc{\C}{{\Bbb C}}
		\nc{\cc}{{\Bbb C}}
\nc{\Cs}{{\Bbb C^*}}
		\nc{\cs}{{\Bbb C^*}}
		\nc{\ccs}{{\Bbb C^*}}
\nc{\D}{{\Bbb D}}
\nc{\E}{{\Bbb E}}
\nc{\F}{{\Bbb F}}
\nc{\G}{{\Bbb G}}
	\nc{\hH}{{\Bbb H}}
\nc{\I}{{\Bbb I}}
\nc{\J}{{\Bbb J}}
\nc{\K}{{\Bbb K}}
	\nc{\lL}{{\Bbb L}}
\nc{\M}{{\Bbb M}}
\nc{\N}{{\Bbb N}}
	\nc{\oO}{{\Bbb O}}
	\nc{\pP}{{\Bbb P}}      
\nc{\Q}{{\Bbb Q}}
\nc{\R}{{\Bbb R}}
	\nc{\sS}{{\Bbb S}}
\nc{\T}{{\Bbb T}}
\nc{\U}{{\Bbb U}}
\nc{\V}{{\Bbb V}}
\nc{\W}{{\Bbb W}}
\nc{\Z}{{\Bbb Z}}
\nc{\X}{{\Bbb X}}
\nc{\Y}{{\Bbb Y}}

\let\P\pP

\nc{\fA}{{\frak A}}
\nc{\fB}{{\frak B}}
\nc{\fC}{{\frak C}}
\nc{\fD}{{\frak D}}
\nc{\fE}{{\frak E}}
\nc{\fF}{{\frak F}}
\nc{\fG}{{\frak G}}
\nc{\fH}{{\frak H}}
\nc{\fI}{{\frak I}}
\nc{\fJ}{{\frak J}}
\nc{\fK}{{\frak K}}
\nc{\fL}{{\frak L}}
\nc{\fM}{{\frak M}}
\nc{\fN}{{\frak N}}
\nc{\fO}{{\frak O}}
\nc{\fP}{{\frak P}}
\nc{\fQ}{{\frak Q}}
\nc{\fR}{{\frak R}}
\nc{\fS}{{\frak S}}
\nc{\fT}{{\frak T}}
\nc{\fU}{{\frak U}}
\nc{\fV}{{\frak V}}
\nc{\fW}{{\frak W}}
\nc{\fZ}{{\frak Z}}
\nc{\fX}{{\frak X}}
\nc{\fY}{{\frak Y}}
\nc{\fa}{{\frak a}}
\nc{\fb}{{\frak b}}
\nc{\fc}{{\frak c}}
\nc{\fe}{{\frak e}}
\nc{\ff}{{\frak f}}
\nc{\fgg}{{\frak g}}
\nc{\fh}{{\frak h}}
\nc{\fiI}{{\frak i}}  
	\nc{\ffi}{{\frak i}}  
\nc{\fj}{{\frak j}}
\nc{\fk}{{\frak k}}
\nc{\fm}{{\frak m}}
\nc{\fn}{{\frak n}}
\nc{\fo}{{\frak o}}
\nc{\fp}{{\frak p}}
\nc{\fq}{{\frak q}}
\nc{\fr}{{\frak r}}
\nc{\fs}{{\frak s}}
\nc{\ft}{{\frak t}}
\nc{\fu}{{\frak u}}
\nc{\fv}{{\frak v}}
\nc{\fw}{{\frak w}}
\nc{\fz}{{\frak z}}
\nc{\fx}{{\frak x}}
\nc{\fy}{{\frak y}}

\nc{\al}{{\alpha }}
\nc{\be}{{\beta }}
\nc{\ga}{{\gamma }}
\nc{\de}{{\delta }}
\nc{\ep}{{\varepsilon }}
\nc{\vap}{{\epsilon }}

\nc{\ze}{{\zeta }}
\nc{\et}{{\eta }}
\rc{\th}{{\theta }}
\nc{\vth}{{\vartheta }}

\nc{\io}{{\iota }}
\nc{\ka}{{\kappa }}
\nc{\la}{{\lambda }}
\nc{\vpi}{{	\varpi		}}
\nc{\vrho}{{	\varrho		}}
\nc{\si}{{	\sigma 		}}
\nc{\ups}{{	\upsilon 	}}
\nc{\vphi}{{	\varphi 	}}
\nc{\om}{{	\omega 		}}

\nc{\Ga}{{\Gamma }}
\nc{\De}{{\Delta }}
\nc{\nab}{{\nabla}}
\nc{\Th}{{\Theta }}
\nc{\La}{{\Lambda }}
\nc{\Si}{{\Sigma }}
\nc{\Ups}{{\Upsilon }}
\nc{\Om}{{\Omega }}

\nc{\Aa}{{	\text{A}	}}
\nc{\Bb}{{	\text{B}	}}
\nc{\Cc}{{	\text{C}	}}
\nc{\Dd}{{	\text{D}	}}
\nc{\Ee}{{	\text{E}	}}
\nc{\Ff}{{	\text{F}	}}
\nc{\Gg}{{	\text{G}	}}
\nc{\Hh}{{	\text{H}	}}
\nc{\Ii}{{	\text{I}	}}
\nc{\Jj}{{	\text{J}	}}
\nc{\Kk}{{	\text{K}	}}
\nc{\Ll}{{	\text{L}	}}
\nc{\Mm}{{	\text{M}	}}
\nc{\Nn}{{	\text{N}	}}
\nc{\Oo}{{	\text{O}	}}
\nc{\Pp}{{	\text{P}	}}
\nc{\Qq}{{	\text{Q}	}}
\nc{\Rr}{{	\text{R}	}}
\nc{\Ss}{{	\text{S}	}}
\nc{\Tt}{{	\text{T}	}}
\nc{\Uu}{{	\text{U}	}}
\nc{\Vv}{{	\text{V}	}}
\nc{\Ww}{{	\text{W}	}}
\nc{\Zz}{{	\text{Z}	}}
\nc{\Xx}{{	\text{X}	}}
\nc{\Yy}{{	\text{Y}	}}

\nc{\bGa}{{	\bbb{\Ga}	}}
\nc{\bA}{{	\bbb{A}		}}
\nc{\bB}{{	\bbb{B}		}}
\nc{\bC}{{	\bbb{C}		}}
\nc{\bD}{{	\bbb{D}		}}
\nc{\bE}{{	\bbb{E}	}}
\nc{\bF}{{	\bbb{F}	}}
\nc{\bG}{{	\bbb{G}	}}
\nc{\bH}{{	\bbb{H}	}}
\nc{\bI}{{	\bbb{I}	}}
\nc{\bJ}{{	\bbb{J}	}}
\nc{\bK}{{	\bbb{K}	}}
\nc{\bL}{{	\bbb{L}	}}
\nc{\bM}{{	\bbb{M}	}}
\nc{\bN}{{	\bbb{N}	}}
\nc{\bO}{{	\bbb{O}	}}
\nc{\bP}{{	\bbb{P}	}}
\nc{\bQ}{{	\bbb{Q}	}}
\nc{\bR}{{	\bbb{R}	}}
\nc{\bS}{{	\bbb{S}	}}
\nc{\bT}{{	\bbb{T}	}}
\nc{\bU}{{	\bbb{U}	}}
\nc{\bV}{{	\bbb{V}	}}
\nc{\bW}{{	\bbb{W}	}}
\nc{\bX}{{	\bbb{X}	}}
\nc{\bY}{{	\bbb{Y}	}}
\nc{\bZ}{{	\bbb{Z}	}}
\nc{\ba}{{	\bbb{a}	}}
			\nc{\bbbb}{{	\bbb{b}	}}
\nc{\bc}{{	\bbb{c}	}}
\nc{\bd}{{	\bbb{d}	}}
			\nc{\bbe}{{	\bbb{e}	}}
			\nc{\bbf}{{	\bbb{f}	}}
\nc{\bg}{{	\bbb{g}	}}
\nc{\bh}{{	\bbb{h}	}}
			\nc{\bbi}{{	\bbb{i}	}}
\nc{\bj}{{	\bbb{j}	}}
			\nc{\bbk}{{	\bbb{k}	}}
\nc{\bl}{{	\bbb{l}	}}
\nc{\bm}{{	\bbb{m}	}}
\nc{\bn}{{	\bbb{n}	}}
\nc{\bo}{{	\bbb{o}	}}
\nc{\bp}{{	\bbb{p}	}}
\nc{\bq}{{	\bbb{q}	}}
\nc{\br}{{	\bbb{r}	}}
\nc{\bs}{{	\bbb{s}	}}
\nc{\bt}{{	\bbb{t}	}}
			\nc{\bbbu}{{	\bbb{u}	}}
\nc{\bv}{{	\bbb{v}	}}
\nc{\bw}{{	\bbb{w}	}}
\nc{\bxx}{{	\bbb{x}	}}
\nc{\by}{{	\bbb{y}	}}
\nc{\bz}{{	\bbb{z}	}}

\nc{\sA}{{\mathsf A}}
\nc{\sB}{{\mathsf B}}
\nc{\sC}{{\mathsf C}}
\nc{\sD}{{\mathsf D}}
\nc{\sE}{{\mathsf E}}
\nc{\sF}{{\mathsf F}}
\nc{\sG}{{\mathsf G}}
\nc{\sH}{{\mathsf H}}
\nc{\sI}{{\mathsf I}}
\nc{\sJ}{{\mathsf J}}
\nc{\sK}{{\mathsf K}}
\nc{\sL}{{\mathsf L}}
\nc{\sM}{{\mathsf M}}
\nc{\sN}{{\mathsf N}}
\nc{\sO}{{\mathsf O}}
\nc{\sP}{{\mathsf P}}
\nc{\sQ}{{\mathsf Q}}
\nc{\sR}{{\mathsf R}}
\rc{\sS}{{\mathsf S}}
\nc{\sT}{{\mathsf T}}
\nc{\sU}{{\mathsf U}}
\nc{\sV}{{\mathsf V}}
\nc{\sW}{{\mathsf W}}
\nc{\sX}{{\mathsf X}}
\nc{\sY}{{\mathsf Y}}
\nc{\sZ}{{\mathsf R}}
\nc{\sa}{{\mathsf a}}
\rc{\sb}{{\mathsf b}}
\rc{\sc}{{\mathsf c}}
\nc{\sd}{{\mathsf d}}
\nc{\sg}{{\mathsf g}}
\nc{\sh}{{\mathsf h}}
\nc{\sj}{{\mathsf j}}
\nc{\sk}{{\mathsf k}}
\nc{\sn}{{\mathsf n}}
\nc{\so}{{\mathsf o}}
\nc{\sq}{{\mathsf q}}
\nc{\sr}{{\mathsf r}}
\nc{\su}{{\mathsf u}}
\nc{\sv}{{\mathsf v}}
\nc{\sw}{{\mathsf w}}
\nc{\sx}{{\mathsf x}}
\nc{\sy}{{\mathsf y}}
\nc{\sz}{{\mathsf z}}

\nc{\toc}{{ 	\small{\tableofcontents} }}
\nc{\addl}{	\addcontentsline{toc}{subsection}	}

\def\sheafEnd{\mathcal{E} \hspace{-1pt} \mathit{nd}}
\def\sheafHom{\mathcal{H} \hspace{-1pt} \mathit{om}}

\newcommand{\bk}{\Bbbk}

\newcommand{\calA}{\mathcal{A}}
\newcommand{\calB}{\mathcal{B}}
\newcommand{\calC}{\mathcal{C}}
\newcommand{\calD}{\mathcal{D}}
\newcommand{\calE}{\mathcal{E}}
\newcommand{\calF}{\mathcal{F}}
\newcommand{\calG}{\mathcal{G}}
\newcommand{\calH}{\mathcal{H}}

\newcommand{\calK}{\mathcal{K}}
\newcommand{\calL}{\mathcal{L}}
\newcommand{\calM}{\mathcal{M}}
\newcommand{\calN}{\mathcal{N}}
\newcommand{\calO}{\mathcal{O}}
\newcommand{\calP}{\mathcal{P}}
\newcommand{\calQ}{\mathcal{Q}}
\newcommand{\calR}{\mathcal{R}}
\newcommand{\calS}{\mathcal{S}}
\newcommand{\calT}{\mathcal{T}}

\newcommand{\calV}{\mathcal{V}}
\newcommand{\calW}{\mathcal{W}}
\newcommand{\calX}{\mathcal{X}}
\newcommand{\calY}{\mathcal{Y}}

\newcommand{\frakg}{\mathfrak{g}}

\newcommand{\frakk}{\mathfrak{k}}

\newcommand{\lotimes}{{\stackrel{_L}{\otimes}}}
\newcommand{\rcap}{{\stackrel{_R}{\cap}}}
\newcommand{\Gm}{{\mathbb{G}}_{\mathbf{m}}}
\newcommand{\Hom}{{\rm Hom}}

\newcommand{\Id}{{\rm Id}}

\newcommand{\Ima}{{\rm Im}}
\newcommand{\Coh}{{\rm Coh}}
\newcommand{\QCoh}{{\rm QCoh}}

\newcommand{\qis}{{\rm qis}}

\newcommand{\fg}{{\rm fg}}
\newcommand{\gr}{{\rm gr}}

\newcommand{\rmi}{\rm{(i)}}
\newcommand{\rmii}{\rm{(ii)}}

\newcommand{\Sym}{{\rm Sym}}

\begin{document}

\begin{abstract}

In this paper we construct, for $F_1$ and $F_2$ subbundles of a vector bundle $E$, a ``Koszul duality'' equivalence between derived categories of $\Gm$-equivariant coherent (dg-)sheaves on the derived intersection $F_1 \rcap_E F_2$, and the corresponding derived intersection $F_1^\perp \rcap_{E^*} F_2^\perp$. We also propose applications to Hecke algebras.

\end{abstract}

\title{Linear Koszul Duality}

\author{Ivan Mirkovi\'c}
\email{mirkovic@math.umass.edu}
\address{University of Massachusetts, Amherst, MA.}

\author{Simon Riche}
\email{riche@math.jussieu.fr}
\address{Universit{\'e} Pierre et Marie Curie, Institut de Math{\'e}matiques de
  Jussieu (UMR 7586 du CNRS), {\'E}quipe d'Analyse Alg{\'e}brique,
175, rue du Cheva\-leret, 75013 Paris, France.}

\maketitle

\section*{Introduction}

\subsection{}

Koszul duality is an algebraic formalism of Fourier transform 
which is 
often deep 
and mysterious in applications. 
For instance, Bezrukavnikov has noticed that
it exchanges monodromy and the Chern class --
the same as mirror duality, while
the work of Beilinson, Ginzburg and Soergel (\cite{BGS}) has
made Koszul duality an essential  ingredient of Representation Theory.

The case of {\em linear Koszul duality} studied here
has a simple geometric content which appears in a number of applications.
For two
vector subbundles 
$F_1, F_2$ 
of a vector bundle $E$ (over a noetherian, integral, separated, regular base scheme),  linear Koszul duality
provides a (contravariant) equivalence of derived categories of $\Gm$-equivariant coherent
sheaves on
the differential graded scheme
\[F_1 \, \rcap_E \, F_2\]
obtained as derived intersection of subbundles inside a vector bundle,
and the corresponding object 
\[F_1^\perp \, \rcap_{E^*} \, F_2^\perp\]
inside the dual vector bundle.

The origin of the linear duality observation
is
Kashiwara's  isomorphism of Borel-Moore homology groups 
\[ \calH_*(F_1 \, \cap_E \, F_2)
\, \cong \,
\calH_*(F_1^\perp \, \cap_{E^*} \, F_2^\perp) \]
given by a Fourier transform for constructible sheaves.
The  
Iwahori-Matsumoto involution
for {\em graded} affine Hecke algebras
has been realized as 
Kashiwara's  Fourier isomorphism
in equivariant Borel-Moore homology
(\cite{EMFou}).
The standard affine Hecke algebras
have analogous realization
in K-theory (the K-homology)
and this suggested that Kashiwara's isomorphism
lifts to K-homology,  but natural isomorphisms
of K-homology groups should come from equivalences of triangulated 
categories of coherent sheaves.

\subsection{} Let us describe the content of this paper. We start in section 1 with generalities on sheaves on dg-schemes. 
In section 2 we construct the relevant Koszul type complexes,
in section 3 we prove the equivalence of categories,
and in section 4 we give the geometric interpretation of this duality.

The idea is that the statement is a particular case of the standard 
Koszul duality in the generality of dg-vector bundles. However, 
because of convergence problems for spectral sequences,
we are able to make sense of this duality only for the
dg-vector bundles with at most 2 non-zero terms. More precisely, our Koszul duality functors are defined in a way similar to those of \cite{GKM}, except for two important differences. First, as explained above, we replace the vector space by a complex of vector bundles. We obtain two ``Koszul dual'' ($\Gm$-equivariant) sheaves of dg-algebras $\calS$ and $\calT$ (which are essentially symmetric algebras of dual complexes of vector bundles.) Then, we modify the functors so that they become \emph{contravariant} and \emph{symmetric}. Indeed, the direct generalization of the constructions of \cite{GKM} would lead us to consider covariant functors of the form \[ \left\{ \begin{array}{ccc} \calS\text{-dg-modules} & \to & \calT\text{-dg-modules} \\ \calM & \mapsto & \calT^{\vee} \otimes_{\calO_X} \calM \end{array} \right. \quad \text{and} \quad \left\{ \begin{array}{ccc} \calT\text{-dg-modules} & \to & \calS\text{-dg-modules} \\ \calN & \mapsto & \calS \otimes_{\calO_X} \calN. \end{array} \right. \] (Here there are some differentials involved, and we have to work with derived functors and derived categories; we do not consider these details in this introduction.) These functors are not well-behaved in general, however, and they are obviously not symmetric for the exchange of $\calS$ and $\calT$. Instead, we consider contravariant functors of the form \[ \left\{ \begin{array}{ccc} \calS\text{-dg-modules} & \to & \calT\text{-dg-modules} \\ \calM & \mapsto & \calT \otimes_{\calO_X} \calM^{\vee} \end{array} \right. \quad \text{and} \quad \left\{ \begin{array}{ccc} \calT\text{-dg-modules} & \to & \calS\text{-dg-modules} \\ \calN & \mapsto & \calS \otimes_{\calO_X} \calN^{\vee}. \end{array} \right. \] The precise definition of these functors is given in subsection \ref{ss:functors}. We use them in subsections \ref{ss:Koszul1} and \ref{ss:Koszul2} to construct some ``generalized Koszul complexes''.

In section 3, we prove that these functors descend to some derived categories, and that they induce equivalences (see subsection \ref{ss:equivalences}). Although the proof is a little technical, its basic idea is very simple: we check (using several spectral sequences) that the composition of these two functors (in any order) is the tensor product with a (generalized) Koszul complex, whose cohomology is trivial. In subsection \ref{ss:restrictionfiniteness} we check that these equivalences respect some finiteness conditions.

Finally, in section 4 we explain the geometric content of these equivalences, \emph{i.e.} we prove that they induce equivalences of categories between $\Gm$-equivariant coherent dg-sheaves on the dg-schemes $F_1 \rcap_E F_2$ and $F_1^{\bot} \rcap_{E^*} F_2^{\bot}$, for subbundles $F_1$ and $F_2$ of a vector bundle $E$.

\subsection{} If we were only interested in characteristic zero, we could have identified the dual of the exterior algebra of a vector bundle with the exterior algebra of the dual vector bundle. Then, for example, the Koszul resolution of the trivial module of the symmetric algebra of a vector space $V$ becomes the symmetric algebra
of the acyclic complex $V \xrightarrow{\Id} V$ (where the first term is in degree $-1$, and the second one in degree $0$). This could have simplified some parts of our constructions. However, in positive characteristic, such an identification is not obvious. Hence we have to pay attention to duals of exterior algebras. In particular, following \cite{BGS}, we rather consider the Koszul resolution above as the tensor product $(\Lambda(V^*))^* \otimes \mathrm{S}(V)$, endowed with a certain differential induced by the natural element in $V^* \otimes V \cong \mathrm{End}(V)$ (see subsection \ref{ss:Koszulcomplexes} for details).  

\subsection{} As explained above, our study involves some derived algebraic geometry. There is a well-developed theory of derived schemes, due to Lurie and Toen, in which derived intersections (and, more generally, derived fibered products) are defined. We have chosen, however, to use the theory of dg-schemes due to Ciocan-Fontanine and Kapranov (see \cite{CK}), which is much more elementary and concrete, and sufficient for our purposes (see also \cite{RKos} and \cite{MR} for other applications of this theory, in representation theoretic contexts). For this purpose, we generalize in section 1 a few well-known facts from the theory of dg-algebras and dg-modules (see \cite{BLEqu}).

\subsection{} A similar geometric interpretation of Koszul duality
has been applied by the second author in \cite{RKos}, in a particular case, to study representations of the Lie algebra
of a connected, simply-connected, semi-simple algebraic group in positive 
characteristic. Let us mention however that the categories and functors
considered here are different from the ones considered in \cite{RKos}. In 
particular, the equivalence of \cite{RKos} is \emph{covariant}, whereas the 
equivalence constructed here is \emph{contravariant}.

\subsection{} In a sequel we will show that
the linear Koszul duality in K-homology is indeed a quantization of Kashiwara's Fourier isomorphism -- the two are related
by the Chern character.
We will also verify that the
linear Koszul duality in equivariant K-homology
gives a geometric realization
of the Iwahori-Matsumoto involution on
(extended) affine Hecke algebras (see \cite{MR}).
This concerns one
typical use of linear Koszul duality.
Consider a partial flag variety
$\calP$ of a group $G$ (either a reductive algebraic group in very good characteristic or a loop group\footnote{Let us point out that the application to loop groups would require an extension of our constructions to the case of infinite dimensional varieties, or ind-schemes, which is not proved here.}),
and
a subgroup $K$
that acts  on $\calP$ with countably many orbits.
Let $\frakg,\frakk$ be the Lie algebras, choose
$E$ to be the trivial bundle $\calP \times \frakg^*$, 
$F_1$ the cotangent subbundle
$T^*\calP$
and $F_2=\calP \times \frakk^\perp$.
Now
$F_1 \, \rcap_E \, F_2$ is a differential graded version of
the Lagrangian $\Lambda_K \subset T^*\calP$, the union of
all conormals to K-orbits in $\calP$,
and $F_1^\perp \, \rcap_{E^*} \, F_2^\perp$ is the
stabilizer dg-scheme for the action of the Lie algebra
$\frakk$ on $\calP$.
If $K$ is the Borel subgroup then 
$F_1^\perp \, \rcap_{E^*} \, F_2^\perp$
is homotopic
to
$F_1 \, \rcap_E \, F_2$
and linear Koszul duality provides
an involution on
the K-group  of
equivariant coherent sheaves on $\Lambda_K$.

Let us conclude by proposing some further applications
of linear  Koszul duality.
The above application to
Iwahori-Matsumoto involutions
should extend to its generalization, the
Aubert involution
on irreducible representations of $p$-adic groups (\cite{AUBDua}).
Linear duality should be an ingredient in
a  geometric realization (proposed in \cite{BFM})
of the
{\em Cherednik Fourier transform}
(essentially an involution on the Cherednik Hecke algebra),
in the
Grojnowski-Garland realization of
Cherednik Hecke algebras as equivariant K-groups
of Steinberg varieties for  affine flag varieties (see \cite{GG},
\cite{VASInd}).
The appearence of linear Koszul duality 
for conormals to Bruhat cells
should also be
a classical limit of the Beilinson-Ginzburg-Soergel Koszul duality
for the mixed category $\calO$ (\cite{BGS}), as 
mixed Hodge modules  come with a deformation
(by Hodge filtration),
to a coherent sheaf
on the  characteristic variety.

\subsection{Acknowledgements}

We are very grateful to Leonid Positselskii, in particular
the essential idea to use contravariant Koszul duality is from
his unpublished lectures at IAS.
We are also grateful to Roman Bezrukavnikov
for his enthusiasm for linear Koszul duality
(which he also named), and to the Institute for Advanced Study
in Princeton for hospitality
and excellent working
environment. Finally the second author thanks Patrick Polo for his helpful remarks.

\section{Generalities on sheaves of dg-algebras and dg-schemes}

In this section $X$ is any noetherian scheme satisfying the following assumption\footnote{See \emph{e.g.} the remarks before \cite[Lemma 2.3.4]{CK} for comments on this assumption.}: $$\begin{array}{c} \text{for any coherent sheaf } \calF \text{ on } X, \text{ there exists a locally free } \\ \text{sheaf of finite rank } \calE \text{ and a surjection } \calE \twoheadrightarrow \calF. \end{array} \leqno(*)$$ We introduce basic definitions concerning dg-schemes and quasi-coherent dg-sheaves, main\-ly following \cite{CK} and \cite{RKos}.

\subsection{Definitions}

Recall the definitions of sheaves of $\calO_X$-dg-algebras and dg-modules given in \cite[1.1]{RKos}.

\begin{defin}

A dg-scheme is a pair $\mathbf{X}=(X,\calA)$ where $X$ is a noetherian scheme satisfying $(*)$, and $\calA$ is a non-positively graded, graded-commutative $\calO_X$-dg-algebra such that $\calA^i$ is a quasi-coherent $\calO_X$-module for any $i \in \mathbb{Z}_{\leq 0}$.

\end{defin}

\begin{defin}

Let $\mathbf{X}=(X,\calA)$ be a dg-scheme.

$\rmi$ A quasi-coherent dg-sheaf $\calF$ on $\mathbf{X}$ is an $\calA$-dg-module such that $\calF^i$ is a quasi-coherent $\calO_X$-module for any $i \in \mathbb{Z}$.

$\rmii$ A coherent dg-sheaf $\calF$ on $\mathbf{X}$ is a quasi-coherent dg-sheaf whose cohomology $H(\calF)$ is a locally finitely generated sheaf of $H(\calA)$-modules.

\end{defin}

We denote by $\calC(\mathbf{X})$, or $\calC(X,\calA)$, the category of quasi-coherent dg-sheaves on the dg-scheme $\mathbf{X}$, and by $\calD(\mathbf{X})$, or $\calD(X,\calA)$, the associated derived category (\emph{i.e.} the localization of the homotopy category of $\calC(\mathbf{X})$ with respect to quasi-isomorphisms).

Similarly, we denote by $\calC^{{\rm c}}(\mathbf{X})$ or $\calC^{{\rm c}}(X,\calA)$, $\calD^{{\rm c}}(\mathbf{X})$ or $\calD^{{\rm c}}(X,\calA)$, the full subcategories whose objects are the coherent dg-sheaves.

If $\mathbf{X}$ is an ordinary scheme, \emph{i.e.} if $\calA=\calO_X$, then we have equivalences $$\calD(\mathbf{X}) \cong \calD \QCoh(X), \quad \calD^{{\rm c}}(\mathbf{X}) \cong \calD^b \Coh(X).$$

Let us stress that these definitions and notation are \emph{different} from the ones used in \cite{RKos} (in $\mathit{loc.} \, \mathit{cit.}$, we only require the cohomology of $\calF$ to be quasi-coherent). This definition will be more suited to our purposes here. Moreover, these two definitions coincide under reasonable assumptions. For the categories of coherent dg-sheaves in all the cases we consider here, this can be deduced from \cite[3.3.4]{RKos}.

\subsection{K-flat resolutions}

Let us fix a dg-scheme $\mathbf{X}=(X,\calA)$. If $\calF$ and $\calG$ are $\calA$-dg-modules, we define as usual the tensor product $\calF \otimes_{\calA} \calG$ (see \cite[1.2]{RKos}). It has a natural structure of an $\calA$-dg-module (here $\calA$ is graded-commutative, hence we do not have to distinguish between left and right dg-modules).

Recall the definition of a K-flat dg-module (see \cite{SPARes}):

\begin{defin}

An $\calA$-dg-module $\calF$ is said to be \emph{K-flat} if for
every $\calA$-dg-module $\calG$ such that $H(\calG)=0$,
we have $H(\calG \otimes_{\calA} \calF)=0$.

\end{defin}

Using \cite[3.4, 5.4.(c)]{SPARes} and assumption $(*)$, one easily proves the following lemma.

\begin{lem}

Let $\calF$ be a quasi-coherent $\calO_X$-dg-module. There exist a quasi-coherent, K-flat $\calO_X$-dg-module $\calP$ and a surjective quasi-isomor\-phism $\calP \xrightarrow{\qis} \calF$.

\end{lem}

Then, using the induction functor $\calF \mapsto \calA \otimes_{\calO_X} \calF$, the following proposition can be proved exactly as in \cite[1.3.5]{RKos}.

\begin{prop} \label{prop:Kflatsqcoh}

Let $\calF$ be a quasi-coherent dg-sheaf on $\mathbf{X}$. There exist a quasi-coherent dg-sheaf $\calP$ on $\mathbf{X}$, K-flat as an $\calA$-dg-module, and a quasi-isomorphism $\calP \xrightarrow{\qis} \calF$.

\end{prop}

\subsection{Invariance under quasi-isomorphisms}

In this subsection we prove that the categories $\calD(\mathbf{X})$, $\calD^{{\rm c}}(\mathbf{X})$ depend on $\calA$ only up to quasi-isomorphism.

Let $X$ be a noetherian scheme satisfying $(*)$, and let $\mathbf{X}=(X,\calA)$ and $\mathbf{X}'=(X,\calB)$ be two dg-schemes with base scheme $X$. Let $\phi: \calA \to \calB$ be a morphism of sheaves of $\calO_X$-dg-algebras. There is a natural functor $$\phi^* : \calC(\mathbf{X}') \to \calC(\mathbf{X})$$ (restriction of scalars), which induces a functor $$R\phi^* : \calD(\mathbf{X}') \to \calD(\mathbf{X}).$$ Similarly, there is a natural functor $$\phi_* : \left\{ \begin{array}{ccc} \calC(\mathbf{X}) & \to & \calC(\mathbf{X}') \\ \calF & \mapsto & \calB \otimes_{\calA} \calF \end{array} . \right.$$ 

We refer to \cite{DELCoh} or \cite{KELUse} for generalities on localization of triangulated categories and derived functors (in the sense of Deligne). The following lemma is borrowed from \cite[5.7]{SPARes} (see also \cite[1.3.6]{RKos}), and implies that K-flat $\calA$-dg-modules are split on the left for the functor $\phi_*$. Using Proposition \ref{prop:Kflatsqcoh}, it follows that $\phi_*$ admits a left derived functor $$L\phi_* : \calD(\mathbf{X}) \to \calD(\mathbf{X}').$$

\begin{lem}

Let $\calF$ be an object of $\calC(X,\calA)$ which is acyclic (\emph{i.e.} $H(\calF)=0$) and K-flat as an $\calA$-dg-module. Then $\calB \otimes_{\calA} \calF$ is acyclic.

\end{lem}

The following result is an immediate extension of \cite[10.12.5.1]{BLEqu}.

\begin{prop} \label{prop:qisqcoh}

$\rmi$ Assume $\phi : \calA \to \calB$ is a quasi-isomorphism. Then the functors $L\phi_*$, $R\phi^*$ are quasi-inverse equivalences of categories $$\calD(\mathbf{X}) \cong \calD(\mathbf{X}').$$ 

$\rmii$ These equivalences restrict to equivalences $$\calD^c(\mathbf{X}) \cong \calD^c(\mathbf{X}').$$

\end{prop}

\emph{Proof}: Statement $\rmi$ can be proved as in \cite[10.12.5.1]{BLEqu} or \cite[1.5.6]{RKos}. Then, clearly, for $\calG$ in $\calD(\mathbf{X}')$ we have $\calG \in \calD^{{\rm c}}(\mathbf{X}')$ iff $R\phi^* \calG \in \calD^{{\rm c}}(\mathbf{X})$. Point $\rmii$ follows. \qed

\subsection{Derived intersection} \label{ss:derivedintersection}

Using Proposition \ref{prop:qisqcoh}, one can consider dg-schemes ``up to quasi-isomorphism'', \emph{i.e.} identify the dg-schemes $(X,\calA)$ and $(X,\calB)$ whenever $\calA$ and $\calB$ are quasi-isomorphic.

As a typical example, we define the derived intersection of two closed subschemes. Consider a scheme $X$, and two closed subschemes
$Y$ and $Z$. Let us denote by $i : Y \hookrightarrow X$ and $j : Z \hookrightarrow X$ the
closed embeddings. Consider the sheaf of dg-algebras $i_* \calO_Y
\lotimes_{\calO_X} j_* \calO_Z$ on $X$. It is well defined up to
quasi-isomorphism: if $\calA_Y \to i_* \calO_Y$, respectively $\calA_Z
\to j_* \calO_Z$ are quasi-isomorphisms of non-positively graded,
graded-commutative sheaves of $\calO_X$-dg-algebras\footnote{See \emph{e.g.} \cite[2.6.1]{CK} for a proof of the existence of such resolutions.}, with $\calA_Y$ and
$\calA_Z$ quasi-coherent and K-flat over $\calO_X$, then $i_* \calO_Y
\lotimes_{\calO_X} j_* \calO_Z$ is quasi-isomorphic to $\calA_Y
\otimes_{\calO_X} j_* \calO_Z$, or to $i_* \calO_Y \otimes_{\calO_X}
\calA_Z$, or to $\calA_Y \otimes_{\calO_X} \calA_Z$.

\begin{defin}

The right derived intersection of $Y$ and $Z$ in $X$ is the dg-scheme
$$Y \, \rcap_X \, Z := (X, \ i_* \calO_Y \, \lotimes_{\calO_X} \, j_* \calO_Z),$$
defined up to quasi-isomorphism.

\end{defin}

To be really precise, only the derived categories $\calD(Y \, \rcap_X \, Z)$, $\calD^{{\rm c}}(Y \, \rcap_X \, Z)$ are well defined (up to equivalence). This is all we will use here.

\section{Generalized Koszul complexes} \label{sec:Koszulcomplexes}

In this section we introduce the dg-algebras we are interested in, and define our Koszul complexes.

\subsection{Notation and definitions} \label{ss:notations}

{From} now on $X$ is a noetherian, integral, separated, regular scheme of dimension $d$. Observe that $X$ satisfies condition $(*)$ by \cite[III.Ex.6.8]{HARAG}. We will consider $\Gm$-equivariant dg-algebras on $X$, \emph{i.e.} sheaves of $\calO_X$-algebras $\calA$, endowed with a $\mathbb{Z}^2$-grading $$\calA = \bigoplus_{i,j \in \mathbb{Z}} \calA^i_j$$ and an $\calO_X$-linear differential $d_{\calA} : \calA \to \calA$, of bidegree $(1,0)$, \emph{i.e.} such that $d_{\calA}(\calA^i_j) \subseteq \calA^{i+1}_j$, and satisfying $$d_{\calA}(a \cdot b)=d_{\calA}(a) \cdot b + (-1)^i a \cdot d_{\calA}(b)$$ for $a \in \calA^i_j$, $b \in \calA$. The basic example is $\calO_X$, endowed with the trivial grading (\emph{i.e.} it is concentrated in bidegree $(0,0)$) and the trivial differential.

A $\Gm$-equivariant dg-module over $\calA$ is a sheaf $\calM$ of $\mathbb{Z}^2$-graded $\calA$-modules endowed with a differential $d_{\calM}$ of bidegree $(1,0)$ satisfying $$d_{\calM}(a \cdot m)=d_{\calA}(a) \cdot m + (-1)^i a \cdot d_{\calM}(m)$$ for $a \in \calA^i_j$, $m \in \calM$. 

We will only consider quasi-coherent ($\Gm$-equivariant) $\calO_X$-dg-algebras. If $\calA$ is such a dg-algebra, we denote by $\calC_{\gr}(\calA)$ the category of quasi-coherent $\Gm$-equivariant $\calA$-dg-modules, \emph{i.e.} $\Gm$-equivariant $\calA$-dg-modules $\calM$ such that $\calM^i_j$ is $\calO_X$-quasi-coherent for any indices $i,j$.

If $\calM$ is a $\Gm$-equivariant $\calA$-dg-module, and $m$ is a local section of $\calM^i_j$, we write $|m|=i$. This integer is called the \emph{cohomological degree} of $m$, while $j$ is called its \emph{internal degree}. We can define two shifts in $\calC_{\gr}(\calA)$: $[n]$, shifting the cohomological degree, and $\langle m \rangle$, shifting the internal degree. More precisely we set $$(\calM[n]\langle m \rangle)^i_j=\calM^{i+n}_{j-m}.$$ Beware that in our conventions $\langle 1 \rangle$ is a ``homological'' shift, \emph{i.e.} it shifts the internal degrees to the right. Also, we use the same conventions as in \cite[{\S}10]{BLEqu} or \cite[1.1]{RKos} concerning the shift [1], \emph{i.e.} the differential of $\calM[1]$ is opposite to the differential of $\calM$.

If $\calM$ and $\calN$ are two $\Gm$-equivariant $\calO_X$-dg-modules, there is a natural structure of $\Gm$-equivariant $\calO_X$-dg-module on the tensor product $\calM \otimes_{\calO_X} \calN,$ with differential defined on homogeneous local sections by $$d_{\calM \otimes \calN}(m \otimes n)=d_{\calM}(m) \otimes n + (-1)^{|m|} m \otimes d_{\calN}(n).$$

If $\calM$ is a $\Gm$-equivariant $\calO_X$-dg-module, we define the $\Gm$-equivariant $\calO_X$-dg-module $\calM^{\vee}$ as the graded dual of $\calM$, \emph{i.e.} the dg-module with $(i,j)$-component $$(\calM^{\vee})^i_j:=\sheafHom_{\calO_X}(\calM^{-i}_{-j}, \calO_X),$$ and with differential defined by $d_{\calM^{\vee}}(f)=-(-1)^{|f|} f \circ d_{\calM}$ for $f \in \calM^{\vee}$ homogeneous. If $\calM$ and $\calN$ are two $\Gm$-equivariant $\calO_X$-dg-modules, there is a natural morphism defined (on homogeneous local sections) by \begin{equation} \label{eq:isomdualtensor} \left\{ \begin{array}{ccc} \calM^{\vee} \otimes_{\calO_X} \calN^{\vee} & \to & (\calM \otimes_{\calO_X} \calN)^{\vee} \\ f \otimes g & \mapsto & \bigl( m \otimes n \mapsto (-1)^{|m| \cdot |g|} f(m) \cdot g(n) \bigr) \end{array} \right. ,\end{equation} which is an isomorphism \emph{e.g.} if the homogeneous components of $\calM$, $\calN$ and $\calM \otimes_{\calO_X} \calN$ are locally free of finite rank. If $\calM$ is a $\Gm$-equivariant $\calO_X$-dg-module such that $\calM^i_j$ is locally-free of finite rank for any $i,j$, then there is an isomorphism \begin{equation} \label{eq:isombidual} \left\{ \begin{array}{ccc} \calM & \xrightarrow{\sim} & (\calM^{\vee})^{\vee} \\ m & \mapsto & \bigl( f \mapsto (-1)^{|f| \cdot |m|} f(m) \bigr) \end{array} \right. . \end{equation}

Let us recall the definition of the truncation functors. If $\calM$ is a $\Gm$-equivariant $\calO_X$-dg-module and if $n \in \mathbb{Z}$, we define the $\Gm$-equivariant $\calO_X$-dg-module $\tau_{\geq n}(\calM)$ by $$\tau_{\geq n}(\calM)^i_j:=\left\{ \begin{array}{cl} 0 & \text{if } i<n \\ \calM^n_j / d_{\calM}(\calM^{n-1}_j) & \text{if } i=n \\ \calM^i_j & \text{if } i>n \end{array} \right. ,$$ with the differential induced by $d_{\calM}$. There is a natural morphism $\calM \to \tau_{\geq n}(\calM)$. Similarly, we define the $\Gm$-equivariant $\calO_X$-dg-module $\tau_{\leq n}(\calM)$ by $$\tau_{\leq n}(\calM):={\rm Ker} \bigl( \calM \to \tau_{\geq n+1}(\calM) \bigr).$$ Observe that if $\calA$ is a $\Gm$-equivariant dg-algebra with $\calA^i_j=0$ for $i>0$, and if $\calM$ is a $\Gm$-equivariant $\calA$-dg-module, then $\tau_{\geq n}(\calM)$ and $\tau_{\leq n}(\calM)$ are again $\Gm$-equivariant $\calA$-dg-modules.

If $\calM$ is a $\Gm$-equivariant $\calO_X$-dg-module, we denote by $\Sym(\calM)$ the graded-symmetric algebra of $\calM$ over $\calO_X$ (\emph{i.e.} the quotient of the tensor algebra of $\calM$ by the relations $m \otimes n = (-1)^{|m| \cdot |n|} n \otimes m$), considered as a $\Gm$-equivariant dg-algebra with differential induced by $d_{\calM}$. Similarly, if $\calF$ is any $\calO_X$-module, we denote by ${\rm S}_{\calO_X}(\calF)$, respectively $\Lambda_{\calO_X}(\calF)$, the symmetric algebra of $\calF$, respectively the exterior algebra of $\calF$, \emph{i.e.} the quotient of the tensor algebra of $\calF$ by the relations $m \otimes n = n \otimes m$, respectively $m \otimes n = - n \otimes m$. Neglecting the gradings, ${\rm S}_{\calO_X}(\calF)$, respectively $\Lambda_{\calO_X}(\calF)$, is the algebra $\Sym(\calF)$, where $\calF$ is concentrated in even cohomological degrees, respectively in odd cohomological degrees. For simplicity, sometimes we drop the subscript ``$\calO_X$''. If $i \geq 0$, we denote by $\mathrm{S}^i(\calF)$, respectively $\Lambda^i(\calF)$, the image of $\calF^{\otimes i}$ in $\mathrm{S}(\calF)$, respectively $\Lambda(\calF)$.

Let us consider two locally free sheaves of finite rank $\calV$ and $\calW$ on $X$, and a morphism of sheaves $f: \calV \to \calW$. 
Let $\calV^{\vee}:=\sheafHom_{\calO_X}(\calV, \calO_X)$ and $\calW^{\vee}:=\sheafHom_{\calO_X}(\calW, \calO_X)$ be the dual locally free sheaves, and $f^{\vee} : \calW^{\vee} \to \calV^{\vee}$ be the morphism induced by $f$. Let us consider the $\Gm$-equivariant $\calO_X$-dg-modules (or complexes of graded $\calO_X$-modules) \[ \calX:= \bigl(\cdots \to 0 \to \calV \xrightarrow{f} \calW \to 0 \to \cdots \bigr), \] where $\calV$ is in bidegree $(-1,2)$ and $\calW$ is in bidegree $(0,2)$, and \[ \calY:= \bigl(\cdots \to 0 \to \calW^{\vee} \xrightarrow{-f^{\vee}} \calV^{\vee} \to 0 \to \cdots \bigr), \] where $\calW^{\vee}$ is in bidegree $(-1,-2)$ and $\calV^{\vee}$ is in bidegree $(0,-2)$.

In sections \ref{sec:Koszulcomplexes} and \ref{sec:duality} we will consider the following $\Gm$-equivariant dg-algebras: \begin{eqnarray*} \calT & := & \Sym( \calX ), \\ \calR & := & \Sym( \calY ), \\ \calS & := & \Sym( \calY[-2]). \end{eqnarray*} For example, the generators of $\calT$ are in bidegrees $(-1,2)$ and $(0,2)$, and the generators of $\calS$ are in bidegrees $(1,-2)$ and $(2,-2)$.

If $\calM$ is a $\Gm$-equivariant $\calS$-dg-module, the dual $\calM^{\vee}$ has a natural structure of a $\calS$-dg-module, constructed as follows. The grading and the differential are defined as above, and the $\calS$-action is defined by $$(s \cdot f)(m)=(-1)^{|s| \cdot |f|}f(s \cdot m),$$ for homogeneous local sections $s$ of $\calS$ and $f$ of $\calM^{\vee}$.

If $\calN$ is a $\calT$-dg-module, respectively a $\calR$-dg-module, the same formulas define on $\calN^{\vee}$ a structure of a $\calT$-dg-module, respectively a $\calR$-dg-module.

\subsection{Reminder on the spectral sequence of a double complex}

Let us recall a few facts on the spectral sequence of a double complex. Let $(C^{p,q})_{p,q \in \mathbb{Z}}$ be a double complex (in any abelian category), with differentials $d'$ (of bidegree $(1,0)$) and $d''$ (of bidegree $(0,1)$). Let ${\rm Tot}(C)$ be the total complex of $C$, \emph{i.e.} the complex with $n$-term $${\rm Tot}(C)^n=\bigoplus_{p+q=n} C^{p,q},$$ and with differential $d'+d''$. The following result is proved \emph{e.g.} in \cite[I.4]{GODTop}.

\begin{prop} \label{prop:spectralsequence}

Assume one of the following conditions is satisfied: \begin{enumerate}
\item There exists $N \in \mathbb{Z}$ such that $C^{p,q}=0$ for $p > N$.
\item There exists $N \in \mathbb{Z}$ such that $C^{p,q}=0$ for $q < N$.
\end{enumerate} Then there is a converging spectral sequence $$E_1^{p,q}=H^q(C^{p,*},d'') \Rightarrow H^{p+q}({\rm Tot}(C)).$$

\end{prop}

\subsection{Reminder on Koszul complexes} \label{ss:Koszulcomplexes}

Let $A$ be a commutative ring, and $V$ be a free $A$-module of finite rank. Let $V^{\vee}=\Hom_A(V,A)$ be the dual $A$-module, and consider the natural morphism \[ i: A \to \Hom_A(V,V) \cong V^{\vee} \otimes_A V,\] sending $1_A$ to  $\Id_V$. Let us first consider the bigraded algebras $\Lambda(V[-1] \langle -2 \rangle)$, the exterior algebra of $V$ placed in bidegree $(1,-2)$, and ${\rm S}(V^{\vee} \langle 2 \rangle)$, the symmetric algebra of $V^{\vee}$ placed in bidegree $(0,2)$. The algebra $\Lambda(V[-1] \langle -2 \rangle)$ acts on the dual $(\Lambda(V[-1] \langle -2 \rangle))^{\vee}$ via $$(t \cdot f)(s)=(-1)^{|t| \cdot |f|} f(ts),$$ where $t,s$ are homogeneous elements of $\Lambda(V[-1] \langle -2 \rangle)$, and $f$ is an homogeneous element of the dual $(\Lambda(V[-1] \langle -2 \rangle))^{\vee}$.

Consider the usual Koszul complex \begin{equation} \label{eq:Koszulcomplex1} {\rm Koszul}_1(V):= {\rm S}(V^{\vee} \langle 2 \rangle) \otimes_A (\Lambda(V[-1] \langle -2 \rangle))^{\vee},\end{equation} where the differential is the composition of the morphism $$\left\{ \begin{array}{ccc}  {\rm S}(V^{\vee}) \otimes_A (\Lambda(V))^{\vee} & \to & {\rm S}(V^{\vee}) \otimes_A (\Lambda(V))^{\vee} \\ s \otimes t & \mapsto & (-1)^{|s|} s \otimes t \end{array} \right. $$ followed by the morphism induced by $i$ $${\rm S}(V^{\vee}) \otimes_A (\Lambda(V))^{\vee} \to {\rm S}(V^{\vee}) \otimes_A V^{\vee} \otimes_A V \otimes_A (\Lambda(V))^{\vee}$$ and finally followed by the morphism $${\rm S}(V^{\vee}) \otimes_A V^{\vee} \otimes_A V \otimes_A (\Lambda(V))^{\vee} \to {\rm S}(V^{\vee}) \otimes_A (\Lambda(V))^{\vee}$$ induced by the action of $V^{\vee} \subset {\rm S}(V^{\vee})$ on ${\rm S}(V^{\vee})$ by right multiplication and the action of $V \subset \Lambda(V)$ on $(\Lambda(V))^{\vee}$ described above. It is well-known (see \emph{e.g.} \cite{BGG}, \cite{BGS}) that this complex has cohomology only in degree $0$, and more precisely that $$H({\rm Koszul_1}(V))=A.$$

The complex ${\rm Koszul_1}(V)$ is a bounded complex of projective graded $A$-modules (here we consider $A$ as a graded ring concentrated in degree $0$). We can take its dual \begin{equation} \label{eq:Koszulcomplex2} {\rm Koszul}_2(V):= ({\rm Koszul}_1(V))^{\vee} \cong \Lambda(V[-1] \langle -2 \rangle) \otimes_A ({\rm S}(V^{\vee} \langle 2 \rangle))^{\vee}. \end{equation} Again we have $$H({\rm Koszul}_2(V))=A.$$

Now, let us consider the bigraded algebras $\Lambda(V[1] \langle -2 \rangle)$, with generators in bidegree $(-1,-2)$, and ${\rm S}(V[-2]\langle 2 \rangle)$, with generators in bidegree $(2,2)$. We have a third Koszul complex \begin{equation} \label{eq:Koszulcomplex3} {\rm Koszul}_3(V):={\rm S}(V^{\vee}[-2]\langle 2 \rangle) \otimes_A (\Lambda(V[1] \langle -2 \rangle))^{\vee}, \end{equation} which may be defined as the bigraded module whose component of bidegree $(i,j)$ satisfies $({\rm Koszul}_3(V))^i_j:=({\rm Koszul}_1(V))^{i-j}_j$, and with differential induced by that of ${\rm Koszul}_1(V)$. As above we have $$H({\rm Koszul}_3(V))=A.$$

We can finally play the same game with the complex ${\rm Koszul}_2(V)$ and obtain the complex \begin{equation} \label{eq:Koszulcomplex4} {\rm Koszul}_4(V) \cong \Lambda(V[1] \langle -2 \rangle) \otimes_A ({\rm S}(V^{\vee} [-2] \langle 2 \rangle))^{\vee} \end{equation} defined by $({\rm Koszul}_4(V))^i_j=({\rm Koszul}_2(V))^{i-j}_j$. Again we have $$H({\rm Koszul_4}(V))=A.$$

\subsection{Two functors} \label{ss:functors}

For any quasi-coherent $\Gm$-equivariant dg-algebra $\calA$ we define the category $\calC^{\searrow}_{\gr}(\calA)$ of $\Gm$-equiva\-riant $\calA$-dg-modules $\calM$ such that $\calM^i_j$ is a coherent $\calO_X$-module for any indices $i$, $j$, and such that there exist integers $N_1, \ N_2$ such that $\calM^i_j=0$ for $i \leq N_1$ or $i+j \geq N_2$. Here the symbol $``\searrow "$ indicates the region in the plane with coordinates $(i,j)$ where the components $\calM^i_j$ can be non-zero, as shown in the figure below.

\begin{figure}[htbp]
\centering
\includegraphics[width=7cm]{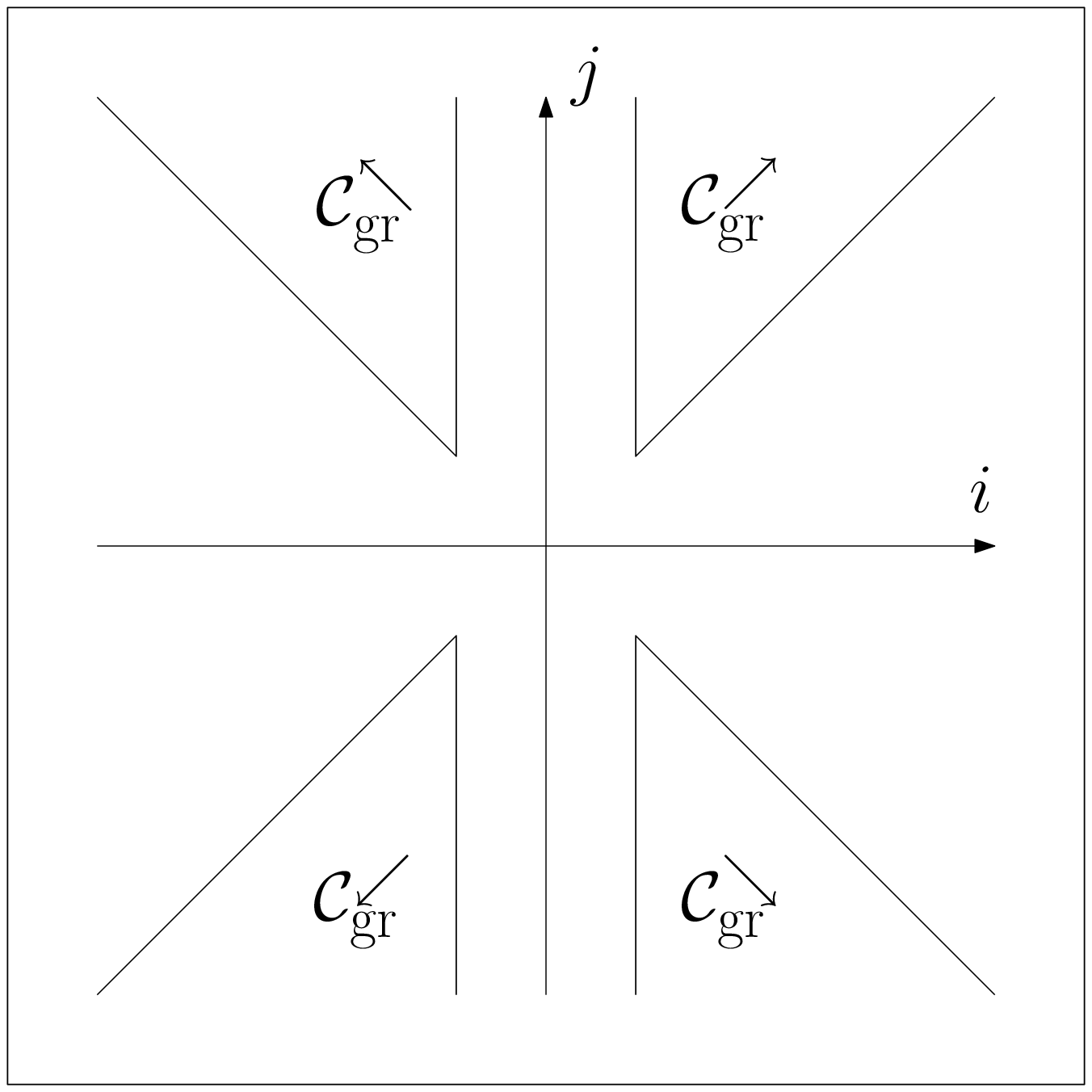}
\end{figure}

Similarly, we define the categories $\calC^{\swarrow}_{\gr}(\calA)$, $\calC^{\nearrow}_{\gr}(\calA)$, $\calC^{\nwarrow}_{\gr}(\calA)$ of $\Gm$-equiva\-riant $\calA$-dg-modules $\calM$ such that the $\calM^i_j$'s are coherent and satisfy the following conditions: $$\begin{array}{cc} \calC^{\swarrow}_{\gr}(\calA) : & \calM^i_j = 0 \text{ if } i \gg 0 \text{ or } i-j \ll 0, \\ \calC^{\nearrow}_{\gr}(\calA) : & \calM^i_j = 0 \text{ if } i \ll 0 \text{ or } i-j \gg 0, \\ \calC^{\nwarrow}_{\gr}(\calA) : & \calM^i_j = 0 \text{ if } i \gg 0 \text{ or } i+j \ll 0. \end{array}$$ 

In this subsection we define two \emph{contravariant} functors $$\mathscr{A}: \calC^{\searrow}_{\gr}(\calS) \to \calC^{\nwarrow}_{\gr}(\calT), \quad \mathscr{B} : \calC^{\nwarrow}_{\gr}(\calT) \to \calC^{\searrow}_{\gr}(\calS).$$ First, let us construct $\mathscr{A}$. If $\calM$ is a $\calS$-dg-module, we have defined in \ref{ss:notations} the $\calS$-dg-module $\calM^{\vee}$. Let $\calM \in \calC_{\gr}^{\searrow}(\calS)$. As a bigraded $\calO_X$-module we set $$\mathscr{A}(\calM)=\calT \otimes_{\calO_X} \calM^{\vee},$$ endowed with a $\calT$-action by left multiplication on the first factor. The differential on $\mathscr{A}(\calM)$ is the sum of four terms. The first one is $d_1:=d_{\calT} \otimes \Id_{\calM^{\vee}}$, and the second one is $d_2:=\Id_{\calT} \otimes d_{\calM^{\vee}}$. Here the tensor product is taken in the graded sense, \emph{i.e.} for homogeneous local sections $t$ and $f$ of $\calT$ and $\calM^{\vee}$ respectively we have $d_2(t \otimes f)=(-1)^{|t|} t \otimes d_{\calM^{\vee}}(f)$. The third and fourth terms are ``Koszul-type'' differentials. Consider first the natural morphism $i: \calO_X \to \sheafEnd_{\calO_X}(\calV) \cong \calV \otimes_{\calO_X} \calV^{\vee}$. Then $d_3$ is the composition of $$\left\{ \begin{array}{ccc} \calT \otimes_{\calO_X} \calM^{\vee} & \to & \calT \otimes_{\calO_X} \calM^{\vee} \\ t \otimes f & \mapsto & (-1)^{|t|} t \otimes f \end{array} \right.$$ followed by the morphism induced by $i$ $$\calT \otimes_{\calO_X} \calM^{\vee} \to \calT \otimes_{\calO_X} \calV \otimes_{\calO_X} \calV^{\vee} \otimes_{\calO_X} \calM^{\vee}$$ and finally followed by the morphism $$\calT \otimes_{\calO_X} \calV \otimes_{\calO_X} \calV^{\vee} \otimes_{\calO_X} \calM^{\vee} \to \calT \otimes_{\calO_X} \calM^{\vee}$$ induced by the right multiplication of $\calV \subset \calT$ on $\calT$, and the left action of $\calV^{\vee} \subset \calS$ on $\calM^{\vee}$. The differential $d_4$ is defined entirely similarly, replacing $\calV$ by $\calW$.

Let us choose a point $x \in X$. Then $\calV_x$, $\calW_x$ are free $\calO_{X,x}$-modules of finite rank. Let $\{v_{\alpha}\}$ be a basis of $\calV_x$, and $\{w_{\beta}\}$ be a basis of $\calW_x$. Let $\{v_{\alpha}^*\}$, $\{w_{\beta}^*\}$ be the dual bases of $(\calV^{\vee})_x$ and $(\calW^{\vee})_x$, respectively. Then the morphism induced by $d_3+d_4$ on $\calT_x \otimes_{\calO_{X,x}} (\calM^{\vee})_x$ can be written \begin{equation} \label{eq:formulad3local} (d_3+d_4)(t \otimes f)=(-1)^{|t|} \bigl( \sum_{\alpha} t v_{\alpha} \otimes v_{\alpha}^* \cdot f + \sum_{\beta} t w_{\beta} \otimes w_{\beta}^* \cdot f \bigr) \end{equation} for homogeneous local sections $t$ of $\calT$ and $f$ of $\calM^{\vee}$.

Using formula \eqref{eq:formulad3local}, one easily checks the relations \begin{equation} \label{eq:relationsd1} (d_1 + d_2)^2=0, \quad (d_3+d_4)^2=0.\end{equation} Further calculations prove the following formula: \begin{equation} \label{eq:relationsd2} (d_1 + d_2) \circ (d_3 + d_4) + (d_3 + d_4) \circ (d_1 + d_2) = 0. \end{equation} It follows from formulas \eqref{eq:relationsd1} and \eqref{eq:relationsd2} that $d_{\mathscr{A}(\calM)}:=d_1 + d_2 + d_3 + d_4$ is indeed a differential. Finally, one easily checks that $\mathscr{A}(\calM)$ is a $\calT$-dg-module, and that it is an object of the category $\calC^{\nwarrow}_{\gr}(\calT)$. Hence the (contravariant) functor $$\mathscr{A} : \calC_{\gr}^{\searrow}(\calS) \to \calC_{\gr}^{\nwarrow}(\calT)$$ is well defined.

Now we define a functor $\mathscr{B}$ in the reverse direction, using similar formulas. Namely if $\calN$ is a $\calT$-dg-module, we have defined above the $\calT$-dg-module $\calN^{\vee}$. If $\calN \in \calC^{\nwarrow}_{\gr}(\calT)$, as a bigraded $\calO_X$-module, we set $$\mathscr{B}(\calN)=\calS \otimes_{\calO_X} \calN^{\vee},$$ and we endow it with the $\calS$-action by left multiplication on the first factor. The differential is again a sum of four terms. The first two are $d_1:=d_{\calS} \otimes \Id_{\calN^{\vee}}$ and $d_2:=\Id_{\calS} \otimes d_{\calN^{\vee}}$. The third one, denoted $d_3$, is defined as above as the composition of $$\left\{ \begin{array}{ccc} \calS \otimes_{\calO_X} \calN^{\vee} & \to & \calS \otimes_{\calO_X} \calN^{\vee} \\ s \otimes g & \mapsto & (-1)^{|s|} s \otimes g \end{array} \right.$$ followed by the morphism induced by $i': \calO_X \to \calV^{\vee} \otimes_{\calO_X} \calV$ $$\calS \otimes_{\calO_X} \calN^{\vee} \to \calS \otimes_{\calO_X} \calV^{\vee} \otimes_{\calO_X} \calV \otimes_{\calO_X} \calN^{\vee}$$ and finally followed by the morphism $$\calS \otimes_{\calO_X} \calV^{\vee} \otimes_{\calO_X} \calV \otimes_{\calO_X} \calN^{\vee} \to \calS \otimes_{\calO_X} \calN^{\vee}$$ induced by the right multiplication of $\calV^{\vee} \subset \calS$ on $\calS$, and the left action of $\calV \subset \calT$ on $\calN^{\vee}$. The differential $d_4$ is defined similarly, replacing $\calV$ by $\calW$. As above, one checks that $d_{\mathscr{B}(\calN)}:=d_1 + d_2 + d_3 + d_4$ is a differential, which turns $\mathscr{B}(\calN)$ into a $\calS$-dg-module, and even an object of $\calC^{\searrow}_{\gr}(\calS)$. For this final claim we use the fact that if $\calS^k_l \neq 0$, then $k+l \leq 0$. As above, this proves that the (contravariant) functor $$\mathscr{B} : \calC^{\nwarrow}_{\gr}(\calT) \to \calC^{\searrow}_{\gr}(\calS)$$ is well defined.

\subsection{First generalized Koszul complex} \label{ss:Koszul1}

Consider the object $$\calK^{(1)}:=\mathscr{B}(\calT) \in \calC_{\gr}^{\searrow}(\calS).$$ It is concentrated in non-negative cohomological degrees, and in non-positive internal degrees.

\begin{lem} \label{lem:Koszul1}

The natural morphism $\calK^{(1)} \to \calO_X$ (projection on the $(0,0)$-component) is a quasi-isomorphism of $\Gm$-equivariant $\calS$-dg-modules.

\end{lem}

\begin{proof} It is sufficient to prove that the localization of this morphism at any $x \in X$ is a quasi-isomorphism. We have isomorphisms \begin{eqnarray*} (\calK^{(1)})_x & \cong & (\calS_x) \otimes_{\calO_{X,x}} \calT^{\vee}_x \\ & \cong & \bigoplus_{i,j,k,l} \Lambda^i (\calW_x^{\vee}) \otimes_{\calO_{X,x}} {\rm S}^j(\calV_x^{\vee}) \otimes_{\calO_{X,x}} (\Lambda^k(\calV_x))^{\vee} \otimes_{\calO_{X,x}} ({\rm S}^l(\calW_x))^{\vee}, \end{eqnarray*} where the symbol $``\vee "$ denotes the dual $\calO_{X,x}$-module, and the term $\Lambda^i (\calW_x^{\vee}) \otimes_{\calO_{X,x}} {\rm S}^j(\calV_x^{\vee}) \otimes_{\calO_{X,x}} (\Lambda^k(\calV_x))^{\vee} \otimes_{\calO_{X,x}} ({\rm S}^l(\calW_x))^{\vee}$ is in cohomological degree $i+2j+k$. The differential on $(\calK^{(1)})_x$ is the sum of four terms: $d_1$, induced by the differential of $\calS_x$; $d_2$, induced by the differential of $\calT_x^{\vee}$; and $d_3$ and $d_4$, the Koszul differentials. The effect of these terms on the indices $i,j,k,l$ may be described as follows: $$d_1 : \left\{ \begin{array}{ccc} i & \mapsto & i-1 \\ j & \mapsto & j+1 \end{array} \right. , \ d_2 : \left\{ \begin{array}{ccc} k & \mapsto & k+1 \\ l & \mapsto & l-1 \end{array} \right. , \ d_3 : \left\{ \begin{array}{ccc} j & \mapsto & j+1 \\ k & \mapsto & k-1 \end{array} \right. , \ d_4 \left\{ \begin{array}{ccc} i & \mapsto & i+1 \\ l & \mapsto & l-1 \end{array} \right. .$$

Disregarding the internal grading, $(\calK^{(1)})_x$ is the total complex of the double complex $(C^{p,q})_{p,q \in \mathbb{Z}}$ whose $(p,q)$-term is $$C^{p,q}:= \bigoplus_{\genfrac{}{}{0pt}{}{p=j+k,}{q=i+j}} \Lambda^i (\calW_x^{\vee}) \otimes_{\calO_{X,x}} {\rm S}^j(\calV_x^{\vee}) \otimes_{\calO_{X,x}} (\Lambda^k(\calV_x))^{\vee} \otimes_{\calO_{X,x}} ({\rm S}^l(\calW_x))^{\vee},$$ and whose differentials are $d'=d_1 + d_2$, $d''=d_3+d_4$. We have $C^{p,q}=0$ if $q<0$, hence by Proposition \ref{prop:spectralsequence} there is a converging spectral sequence $$E_1^{p,q}=H^q(C^{p,*},d'') \Rightarrow H^{p+q}((\calK^{(1)})_x).$$ It follows that, to prove the lemma, we only have to prove that the cohomology of $\calS_x \otimes_{\calO_{X,x}} \calT_x^{\vee}$ with respect to the differential $d_3+d_4$ is $\calO_{X,x}$ in degree $0$, and $0$ in other degrees. But this complex is the tensor product of the Koszul complexes ${\rm Koszul}_3(\calV_x)$ (with the internal grading opposite to that in \eqref{eq:Koszulcomplex3}) and ${\rm Koszul}_2(\calW^{\vee}_x)$ of \eqref{eq:Koszulcomplex2}, both living in non-negative degrees. We have seen that these complexes have cohomology $\calO_{X,x}$, and their components are free (hence flat). The result follows, using K\"{u}nneth formula. \end{proof}

\subsection{Second generalized Koszul complex} \label{ss:Koszul2}

Consider now the object $$\calK^{(2)}:=\mathscr{A}(\calS) \in \calC_{\gr}^{\nwarrow}(\calT).$$ It is concentrated in non-positive cohomological degrees, and in non-negative internal degrees. As in \ref{ss:Koszul1}, we are going to prove:

\begin{lem} \label{lem:Koszul2}

The natural morphism $\calK^{(2)} \to \calO_X$ (projection on the $(0,0)$-component) is a quasi-isomorphism of $\Gm$-equivariant $\calT$-dg-modules.

\end{lem}

\begin{proof} The arguments for this proof are completely similar to those of Lemma \ref{lem:Koszul1}. Here the double complex to consider has $(p,q)$-term $$C^{p,q}:=\bigoplus_{\genfrac{}{}{0pt}{}{p=-i-l,}{q=-k-l}} \Lambda^i (\calV_x) \otimes_{\calO_{X,x}} {\rm S}^j(\calW_x) \otimes_{\calO_{X,x}} (\Lambda^k(\calW^{\vee}_x))^{\vee} \otimes_{\calO_{X,x}} ({\rm S}^l(\calV^{\vee}_x))^{\vee}$$ and differentials $d'=d_1+d_2$, $d''=d_3+d_4$. We have $C^{p,q}=0$ for $p>0$.\end{proof}

\section{Algebraic duality} \label{sec:duality}

In this section we prove our Koszul duality between $\calS$- and $\calT$-dg-modules.

\subsection{Resolutions}

First we need to prove the existence of some resolutions.

\begin{prop} \label{prop:resolutions}

$\rmi$ Let $\calM$ be an object of $\calC_{\gr}^{\searrow}(\calS)$. There exist an object $\calP$ of $\calC_{\gr}^{\searrow}(\calS)$ such that, for all indices $i$ and $j$, $\calP^i_j$ is $\calO_X$-locally free of finite rank, and a quasi-isomorphism of $\calS$-dg-modules $\calP \xrightarrow{\qis} \calM$.

$\rmii$ Let $\calN$ be an object of $\calC_{\gr}^{\nwarrow}(\calT)$. There exist an object $\calQ$ of $\calC^{\nwarrow}_{\gr}(\calT)$ such that, for all indices $i$ and $j$, $\calQ^i_j$ is $\calO_X$-locally free of finite rank, and a quasi-isomorphism of $\calT$-dg-modules $\calQ \xrightarrow{\qis} \calN$.

\end{prop}

\begin{proof} We give a proof only for point $\rmi$. The proof of $\rmii$ is similar\footnote{One could also use the ``regrading trick'' of \ref{ss:regrading} below to show that these two statements are equivalent.}. Let $\calM$ be an object of $\calC^{\searrow}_{\gr}(\calS)$. Let $N_1$ and $N_2$ be integers such that $\calM^i_j=0$ for $i<N_1$ or $i+j>N_2$. First we consider $\calM$ as a $\Gm$-equivariant $\calO_X$-dg-module. Then, for each $j \leq N_2 - N_1$, $\calM_j$ is a complex of coherent $\calO_X$-modules, with non-zero terms only in the interval $[N_1, N_2 - j ]$ (and $\calM_j=0$ otherwise). Using a standard procedure (see \emph{e.g.} \cite[I.4.6]{HARRD} and \cite[III.Ex.6.9]{HARAG}), there exists a complex $\calL_j$ of locally free $\calO_X$-modules of finite rank, with non-zero terms only in the interval $[N_1, N_2 - j ]$, and a \emph{surjective} morphism of $\calO_X$-dg-modules $\calL_j \twoheadrightarrow \calM_j$. Then $\calL:=\bigoplus_j \calL_j$ is an object of $\calC^{\searrow}_{\gr}(\calO_X)$, and there is a surjective morphism of $\Gm$-equivariant $\calO_X$-dg-modules $\calL \twoheadrightarrow \calM$. Then $\calP^{(1)}:=\calS \otimes_{\calO_X} \calL$, endowed with the natural differential and the natural action of $\calS$, is an object of $\calC_{\gr}^{\searrow}(\calS)$, and there is a surjective morphism of $\Gm$-equivariant $\calS$-dg-modules $$\calP^{(1)} \twoheadrightarrow \calM.$$

Taking the kernel of this morphism, and repeating the procedure, we obtain objects $\calP^{(i)}$ ($i=1, \cdots, d$) of $\calC_{\gr}^{\searrow}(\calS)$, (recall that $d=\dim(X)$) whose homogeneous components are locally free of finite rank over $\calO_X$, and an exact sequence of $\calS$-dg-modules $$\calP^{(d)} \to \calP^{(d-1)} \to \cdots \to \calP^{(1)} \to \calM \to 0.$$ We define $\calP^{(d+1)}:={\rm ker}(\calP^{(d)} \to \calP^{(d-1)})$. Then, for any indices $i,j$, the exact sequence $$0 \to (\calP^{(d+1)})^i_j \to \cdots \to (\calP^{(1)})^i_j \to \calM^i_j \to 0$$ is a resolution of the $\calO_X$-coherent sheaf $\calM^i_j$, the terms $(\calP^{(k)})^i_j$ being locally free of finite rank for $k=1, \cdots, d$. It follows that $(\calP^{(d+1)})^i_j$ is also locally free of finite rank over $\calO_X$ (see again \cite[III.Ex.6.9]{HARAG}).

Finally we take $$\calP:={\rm Tot} \bigl( 0 \to \calP^{(d+1)} \to \calP^{(d)} \to \cdots \to \calP^{(1)} \to 0 \bigr).$$ It is naturally an object of $\calC^{\searrow}_{\gr}(\calS)$, and an easy spectral sequence argument shows that the natural morphism $\calP \to \calM$ is a quasi-isomorphism of $\calS$-dg-modules. \end{proof}

\subsection{Derived functors} \label{ss:derivedfunctors}

Let us introduce some notation. If $\calA$ is any quasi-coherent $\Gm$-equivariant dg-algebra, we denote by $\calH^{*}_{\gr}(\calA)$ the homotopy category of the category $\calC_{\gr}^{*}(\calA)$, where $*=\nearrow, \nwarrow, \swarrow, \searrow$. The objects of $\calH^*_{\gr}(\calA)$ are the same as those of $\calC_{\gr}^*(\calA)$, and the morphisms in $\calH_{\gr}^*(\calA)$ are the quotient of the morphisms in $\calC^*_{\gr}(\calA)$ by the homotopy relation. These categories are naturally triangulated. We denote by $\calD^*_{\gr}(\calA)$, the localization of $\calH_{\gr}^*(\calA)$ with respect to quasi-isomorphisms.

As a corollary of Proposition \ref{prop:resolutions}, we obtain the following result.

\begin{cor}

The functors $\mathscr{A}$ and $\mathscr{B}$ admit derived functors (in the sense of Deligne) $$\overline{\mathscr{A}} : \calD^{\searrow}_{\gr}(\calS) \to \calD_{\gr}^{\nwarrow}(\calT),  \qquad \overline{\mathscr{B}} : \calD_{\gr}^{\nwarrow}(\calT) \to \calD_{\gr}^{\searrow}(\calS).$$

\end{cor}

\begin{remark} The functor $\overline{\mathscr{A}}$ is the left derived functor of $\mathscr{A}$ if we consider it as a covariant functor $\calC_{\gr}^{\searrow}(\calS) \to \calC_{\gr}^{\nwarrow}(\calT)^{{\rm opp}}$, or the right derived functor of $\mathscr{A}$ if we consider it as a covariant functor $\calC_{\gr}^{\searrow}(\calS)^{{\rm opp}} \to \calC_{\gr}^{\nwarrow}(\calT)$. \end{remark}

\begin{proof} \emph{Case of the functor $\mathscr{A}$.} To fix notations, in this proof we consider $\mathscr{A}$ as a covariant functor $\calC_{\gr}^{\searrow}(\calS) \to \calC_{\gr}^{\nwarrow}(\calT)^{{\rm opp}}$. To prove that $\mathscr{A}$ admits a left derived functor, it is enough to prove that there are enough objects split on the left\footnote{Recall that an object $\calM$ of $\calC_{\gr}^{\searrow}(\calS)$ is said to be \emph{split on the left} for $\mathscr{A}$ if for any quasi-isomorphism $\calM' \xrightarrow{\qis} \calM$, there exists an object $\calM''$ of $\calC_{\gr}^{\searrow}(\calS)$ and a quasi-isomorphism $\calM'' \xrightarrow{\qis} \calM'$ such that the induced morphism $\mathscr{A}(\calM'') \to \mathscr{A}(\calM)$ is again a quasi-isomorphism.} for $\mathscr{A}$ in the category $\calC^{\searrow}_{\gr}(\calS)$ (see \cite{DELCoh} or \cite{KELUse}). To prove the latter fact, using Proposition \ref{prop:resolutions}$\rmi$, it is enough to prove that if $f: \calP \to \calQ$ is a quasi-isomorphism between two objects of $\calC^{\searrow}_{\gr}(\calS)$ whose homogeneous components are $\calO_X$-locally free of finite rank, then the induced morphism $$\mathscr{A}(f): \mathscr{A}(\calP) \to \mathscr{A}(\calQ)$$ is again a quasi-isomorphism. Taking cones, this amounts to proving that if $\calP$ is an acyclic object of $\calC^{\searrow}_{\gr}(\calS)$ whose homogeneous components are $\calO_X$-locally free of finite rank, then $\mathscr{A}(\calP)$ is again acyclic.

So, let $\calP$ be such a $\Gm$-equivariant $\calS$-dg-module. For each index $j$, the complex of $\calO_X$-modules $\calP_j$ is acyclic, bounded, and all its components are locally free of finite rank. It follows that $\calP^{\vee}$ is also acyclic. Let $x$ be a point of $X$, and let us prove that $\mathscr{A}(\calP)_x$ is acyclic. We use the same notations as in \ref{ss:functors}. In particular, $d_{\mathscr{A}(\calP)}$ is the sum of four terms $d_1$, $d_2$, $d_3$ and $d_4$. We have an isomorphism $$\mathscr{A}(\calP)_x \cong \bigoplus_{i,j,k,l} \Lambda^i(\calV_x) \otimes_{\calO_{X,x}} {\rm S}^j(\calW_x) \otimes_{\calO_{X,x}} (\calP^{\vee}_x)^k_l,$$ where the term $\Lambda^i(\calV_x) \otimes_{\calO_{X,x}} {\rm S}^j(\calW_x) \otimes_{\calO_{X,x}} (\calP^{\vee}_x)^k_l$ is in cohomological degree $k-i$. The effect of the differentials on the indices $i,j,k,l$ may be described as $$d_1 : \left\{ \begin{array}{ccc} i & \mapsto & i-1 \\ j & \mapsto & j+1 \end{array} \right. , \quad d_2 : k \mapsto k+1 , \quad d_3 : \left\{ \begin{array}{ccc} i & \mapsto & i+1 \\k & \mapsto & k+2 \\ l & \mapsto & l-2 \end{array} \right. , \quad d_4: \left\{ \begin{array}{ccc} j & \mapsto & j+1 \\ k & \mapsto & k+1 \\ l & \mapsto & l-2 \end{array} \right. .$$

Hence, disregarding the internal grading, $\mathscr{A}(\calP)_x$ is the total complex of the double complex with $(p,q)$-term $$C^{p,q}:= \bigoplus_{\genfrac{}{}{0pt}{}{p=-i-j-l,}{q=k+l+j}} \Lambda^i(\calV_x) \otimes_{\calO_{X,x}} {\rm S}^j(\calW_x) \otimes_{\calO_{X,x}} (\calP^{\vee}_x)^k_l,$$ with differentials $d'=d_3+d_4$ and $d''=d_1 + d_2$. By definition, $\calP$ is in $\calC^{\searrow}_{\gr}(\calS)$, hence $(\calP^{\vee})^k_l=0$ for $k+l \ll 0$. Hence $C^{p,q}=0$ for $q \ll 0$. By Proposition \ref{prop:spectralsequence}, it follows that there is a converging spectral sequence $$E_1^{p,q}=H^q(C^{p,*},d'') \Rightarrow H^{p+q}(\mathscr{A}(\calP)_x).$$ Hence we can forget about the differentials $d_3$ and $d_4$, \emph{i.e.} it is sufficient to prove that the tensor product of $\calO_{X,x}$-dg-modules $$\calT_x \otimes_{\calO_{X,x}} \calP_x^{\vee}$$ is acyclic. We have seen above that $\calP_x^{\vee}$ is acyclic, and $\calT_x$ is a bounded complex of flat $\calO_{X,x}$-modules. Hence $\calT_x \otimes_{\calO_{X,x}} \calP_x^{\vee}$ is indeed acyclic, which finishes the proof of the existence of the derived functor $$\overline{\mathscr{A}} : \calD^{\searrow}_{\gr}(\calS) \to \calD_{\gr}^{\nwarrow}(\calT).$$

\emph{Case of the functor $\mathscr{B}$.} The proof for the functor $\mathscr{B}$ is very similar. If $\calQ$ is a $\Gm$-equivariant $\calT$-dg-module as in Proposition \ref{prop:resolutions}$\rmii$ which is acyclic, and $x \in X$, then we have $$\mathscr{B}(\calQ)_x = \bigoplus_{i,j,k,l} \Lambda^i(\calW_x^{\vee}) \otimes_{\calO_{X,x}} {\rm S}^j(\calV_x^{\vee}) \otimes_{\calO_{X,x}} (\calQ^{\vee}_x)^k_l,$$ where the term $\Lambda^i(\calW_x^{\vee}) \otimes_{\calO_{X,x}} {\rm S}^j(\calV_x^{\vee}) \otimes_{\calO_{X,x}} (\calQ^{\vee}_x)^k_l$ is in cohomological degree $i+2j+k$. Again $\calQ^{\vee}$ is acyclic, and $d_{\mathscr{B}(\calN)}$ is the sum of four terms $d_1$, $d_2$, $d_3$ and $d_4$, whose effect on the indices $i,j,k,l$ may be described as $$d_1 : \left\{ \begin{array}{ccc} i & \mapsto & i-1 \\ j & \mapsto & j+1 \end{array} \right. , \quad d_2 : k \mapsto k+1 , \quad d_3 : \left\{ \begin{array}{ccc} j & \mapsto & j+1 \\k & \mapsto & k-1 \\ l & \mapsto & l+2 \end{array} \right. , \quad d_4: \left\{ \begin{array}{ccc} i & \mapsto & i+1 \\ l & \mapsto & l+2 \end{array} \right. .$$ Hence, disregarding the internal grading, $\mathscr{B}(\calQ)_x$ is the total complex of the double complex with $(p,q)$-term $$D^{p,q}:=\bigoplus_{\genfrac{}{}{0pt}{}{p=i+j,}{q=k+j}} \Lambda^i(\calW_x^{\vee}) \otimes_{\calO_{X,x}} {\rm S}^j(\calV_x^{\vee}) \otimes_{\calO_{X,x}} (\calQ_x^{\vee})^k_l,$$ and with differentials $d'=d_3+d_4$, $d''=d_1 + d_2$. We know that $(\calQ_x^{\vee})^k_l=0$ if $k \ll 0$, hence $D^{p,q}=0$ for $p \ll 0$. By Proposition \ref{prop:spectralsequence}, it follows that there is a converging spectral sequence $$E_1^{p,q}=H^q(D^{p,*},d'') \Rightarrow H^{p+q}(\mathscr{B}(\calQ)_x).$$ Hence it is sufficient to prove that the tensor product of $\calO_{X,x}$-dg-modules $$\calS_x \otimes_{\calO_{X,x}} \calQ_x^{\vee}$$ is acyclic. 

The ($\Gm$-equivariant) $\calO_{X,x}$-dg-module $\calS_x$ has a finite filtration with subquotients finite numbers of copies of ${\rm S}(\calV_x^{\vee})$. Hence it is enough to prove that ${\rm S}(\calV_x^{\vee}) \otimes_{\calO_{X,x}} \calQ_x^{\vee}$ is acyclic. But ${\rm S}(\calV_x^{\vee})$, as a ($\Gm$-equivariant) $\calO_{X,x}$-dg-module, is a direct sum of flat $\calO_{X,x}$-modules (placed in different degrees), hence the latter fact is clear. \end{proof}

\subsection{Morphisms of functors} \label{ss:morphismfunctors}

In this subsection we construct some morphisms of functors. We will prove in the next subsection that they are isomorphisms, which implies that $\overline{\mathscr{A}}$ and $\overline{\mathscr{B}}$ are equivalences of categories.

\begin{prop} \label{prop:morphismsfunctors}

There exist natural morphisms of functors $$\overline{\mathscr{B}} \circ \overline{\mathscr{A}} \to \Id_{\calD^{\searrow}_{\gr}(\calS)}, \qquad \overline{\mathscr{A}} \circ \overline{\mathscr{B}} \to \Id_{\calD^{\nwarrow}_{\gr}(\calT)}.$$

\end{prop}

\begin{proof} Let us give the details for the first morphism. The construction of the second one is similar. It is sufficient to construct this morphism for any $\calA$-dg-module $\calP$ as in Proposition \ref{prop:resolutions}$\rmi$. In this case $\overline{\mathscr{A}}(\calP)$ is isomorphic to the image of $\mathscr{A}(\calP)$ in the derived category. As $\mathscr{A}(\calP)$ has also $\calO_X$-locally free homogeneous components, $\overline{\mathscr{B}} \circ \overline{\mathscr{A}}(\calP)$ is isomorphic to the image of $\mathscr{B} \circ \mathscr{A}(\calP)$ in the derived category. We will define a morphism in $\calC^{\searrow}_{\gr}(\calS)$ \begin{equation} \label{eq:morphismBA} \mathscr{B} \circ \mathscr{A} (\calP) \to \calP. \end{equation}

First we begin with the following lemma, which can be checked by direct computation, using the isomorphisms \eqref{eq:isomdualtensor} and \eqref{eq:isombidual}.

\begin{lem} \label{lem:dualA}

As a bigraded $\calO_X$-module, $(\mathscr{A}(\calP))^{\vee}$ is naturally isomorphic to $\calT^{\vee} \otimes_{\calO_X} \calP$. Under this isomorphism, locally around a point $x \in X$, with the notation of {\rm \eqref{eq:formulad3local}}, the differential becomes \begin{equation*} d_{(\mathscr{A}(\calP))^{\vee}}(f \otimes p) \ = \ d(f) \otimes p + (-1)^{|f|} f \otimes d(p) \ - \ (-1)^{|f|} \bigl( \sum_{\alpha} f \cdot v_{\alpha} \otimes v_{\alpha}^* \cdot p + \sum_{\beta} f \cdot w_{\beta} \otimes w_{\beta}^* \cdot p \bigr),\end{equation*} where we set $(f \cdot t)(t')=f(t \cdot t')$ for $f \in \calT^{\vee}$ and $t,t' \in \calT$.

\end{lem}

Under the isomorphism of Lemma \ref{lem:dualA}, we have as bigraded $\calO_X$-modules $$ \mathscr{B} \circ \mathscr{A} (\calP) \cong \calS \otimes_{\calO_X} \calT^{\vee} \otimes_{\calO_X} \calP.$$ We define the morphism of bigraded $\calO_X$-modules $$\left\{\ \begin{array}{ccc} \calS \otimes_{\calO_X} \calT^{\vee} \otimes_{\calO_X} \calP & \to & \calP \\ s \otimes f \otimes p & \mapsto & f(1_{\calT}) \cdot s \cdot p \end{array} \right. .$$ This morphism clearly commutes with the $\calS$-actions. Moreover, using Lemma \ref{lem:dualA}, one easily checks that it also commutes with the differentials, hence defines the desired morphism \eqref{eq:morphismBA}.\end{proof}

\subsection{Equivalences} \label{ss:equivalences}

\begin{thm}

The functors $\overline{\mathscr{A}}, \ \overline{\mathscr{B}}$ are equivalences of categories, quasi-inverse to each other.

\end{thm}

\begin{proof} \emph{First step: isomorphism $\overline{\mathscr{B}} \circ \overline{\mathscr{A}} \xrightarrow{\sim} \Id$.} In Proposition \ref{prop:morphismsfunctors}, we have constructed a morphism of functors $\overline{\mathscr{B}} \circ \overline{\mathscr{A}} \to \Id$. In this first step we prove that it is an isomorphism. Let $\calP$ be an object of $\calC^{\searrow}_{\gr}(\calS)$ as in Proposition \ref{prop:resolutions}$\rmi$. We have seen in \ref{ss:morphismfunctors} that $\overline{\mathscr{B}} \circ \overline{\mathscr{A}}(\calP)$ is isomorphic to the image of $\mathscr{B} \circ \mathscr{A}(\calP)$ in the derived category. By Proposition \ref{prop:resolutions}$\rmi$, it is thus enough to prove that the induced morphism $$\phi: \mathscr{B} \circ \mathscr{A}(\calP) \to \calP$$ is a quasi-isomorphism. Let us construct a section (over $\calO_X$) for this morphism. As a bigraded $\calO_X$-module we have $\mathscr{B} \circ \mathscr{A}(\calP) \cong \calS \otimes_{\calO_X} \calT^{\vee} \otimes_{\calO_X} \calP$. Let $\epsilon_{\calT} \in \calT^{\vee}$ be the unit section in $(\calT^{\vee})^0_0=\calO_X$. Now consider the morphism $$\psi: \left\{ \begin{array}{ccc} \calP & \to & \mathscr{B} \circ \mathscr{A}(\calP) \\ p & \mapsto & 1_{\calS} \otimes \epsilon_{\calT} \otimes p \end{array} \right. .$$ One easily checks that it is a morphism of $\Gm$-equivariant $\calO_X$-dg-modules (but of course not of $\calS$-dg-modules), and that $$\phi \circ \psi = \Id_{\calP}.$$ Hence it is enough to prove that $\psi$ is a quasi-isomorphism.

As a bigraded $\calO_X$-module, we have, with the notation of \ref{ss:Koszul1}, $$\mathscr{B} \circ \mathscr{A}(\calP) \cong \calK^{(1)} \otimes_{\calO_X} \calP \cong \bigoplus_{i,j,k,l} (\calK^{(1)})^i_k \otimes_{\calO_X} \calP^j_l,$$ where the term $(\calK^{(1)})^i_k \otimes_{\calO_X} \calP^j_l$ is in cohomological degree $i+j$. Remark that here the non-zero terms occur only when $k$ is even. By Lemma \ref{lem:dualA}, the differential on $\mathscr{B} \circ \mathscr{A}(\calP)$ is the sum of four terms. The first one is $d_1:=d_{\calK^{(1)}} \otimes \Id_{\calP}$. The second one is $d_2:=\Id_{\calK^{(1)}} \otimes d_{\calP}$. The third one is the ``Koszul type" differential coming from the left action of $\calV^{\vee} \subset \calS$ on $\calP$ and the right action of $\calV \subset \calT$ on $\calK^{(1)}$. Finally $d_4$ is the similar ``Koszul-type" differential coming from the actions of $\calW^{\vee}$ and $\calW$. The effect of these differentials on the indices $i,j,k,l$ can be described as follows: $$d_1 : i \to i+1, \quad d_2 : j \mapsto j+1, \quad d_3: \left\{ \begin{array}{ccc} i & \mapsto & i-1 \\ j & \mapsto & j+2 \\ k & \mapsto & k+2 \\ l & \mapsto & l-2 \end{array} \right. , \quad d_4: \left\{ \begin{array}{ccc} j & \mapsto & j+1 \\ k & \mapsto & k+2 \\ l & \mapsto & l-2 \end{array} \right. .$$ Moreover, one easily checks the following relations: $$(d_1 + d_4)^2=0, \qquad (d_2 + d_3)^2=0.$$ Hence, disregarding the internal grading, $\mathscr{B} \circ \mathscr{A}(\calP)$ is the total complex of the double complex with $(p,q)$-term $$C^{p,q}:=\bigoplus_{\genfrac{}{}{0pt}{}{p=j+l+k/2,}{q=i-l-k/2}} (\calK^{(1)})^i_k \otimes_{\calO_X} \calP^j_l,$$ and with differentials $d'=d_2+d_3$ and $d''=d_1+d_4$. We know that $\calP^j_l=0$ for $j+l \gg 0$, and that $(\calK^{(1)})^i_k=0$ if $k > 0$. Hence $C^{p,q}=0$ for $p \gg 0$. It follows, by Proposition \ref{prop:spectralsequence}, that there is a converging spectral sequence $$E_1^{p,q}=H^q(C^{p,*},d'') \Rightarrow H^{p+q}(\mathscr{B} \circ \mathscr{A} (\calP)).$$ Disregarding the internal grading, $\calP$ is also the total complex of a double complex, defined by $$(C')^{p,q}:=\calP^{p+q}_{-q}$$ and the differentials $d'=d_{\calP}$, $d''=0$. Here also $(C')^{p,q}=0$ for $p \gg 0$, hence the corresponding spectral sequence converges. Moreover, $\psi$ is induced by a morphism of double complexes $C' \to C$. It follows that it is enough to prove that the morphism induced by $\psi$ from $\calP$, endowed with the zero differential, to $\calK^{(1)} \otimes_{\calO_X} \calP$, endowed with the differential $d_1+d_4$, is a quasi-isomorphism.

The latter dg-module is again the total complex of the double complex with $(p,q)$-term $$D^{p,q}:=\bigoplus_{k,l} (\calK^{(1)})^q_k \otimes_{\calO_X} \calP^p_l,$$ and differentials $d'=d_4$, $d''=d_1$. And $\calP$ (with the trivial differential) is also the total complex of the double complex defined by $$(D')^{p,q}=\left\{ \begin{array}{cl} \bigoplus_l \calP^p_l & \text{if } q=0, \\[2pt] 0 & \text{otherwise}, \end{array} \right.$$ and with two trivial differentials. Again $\psi$ is induced by a morphism of double complexes, and we have $D^{p,q}=(D')^{p,q}=0$ for $q<0$. We conclude that the associated spectral sequences converge. As $H(\calK_1)=\calO_X$ (see Lemma \ref{lem:Koszul1}) and $\calP$ is a bounded above complex of flat $\calO_X$-modules, we finally conclude that $\psi$ is a quasi-isomorphism. \bigskip

\emph{Second step: isomorphism $\overline{\mathscr{A}} \circ \overline{\mathscr{B}} \xrightarrow{\sim} \Id$}. The proofs in this second step are very similar to those of the first step. By Proposition \ref{prop:morphismsfunctors} there is a natural morphism $\overline{\mathscr{A}} \circ \overline{\mathscr{B}} \to \Id$, and we prove that it is an isomorphism. As above, it is enough to prove that, for $\calQ$ an object of $\calC_{\gr}^{\nwarrow}(\calT)$ as in Proposition \ref{prop:resolutions}$\rmii$, the induced morphism of $\calT$-dg-modules $$\phi: \mathscr{A} \circ \mathscr{B}(\calQ) \to \calQ$$ is a quasi-isomorphism. Also as above one can construct a section $$\psi : \calQ \to \mathscr{A} \circ \mathscr{B}(\calQ)$$ of $\phi$ as a morphism of $\Gm$-equivariant $\calO_X$-dg-modules, and it is enough to prove that $\psi$ is a quasi-isomorphism.

Here we have as bigraded $\calO_X$-modules, with the notation of \ref{ss:Koszul2}, $$\mathscr{A} \circ \mathscr{B}(\calQ) \cong \calK^{(2)} \otimes_{\calO_X} \calQ \cong \bigoplus_{i,j,k,l} (\calK^{(2)})^i_k \otimes_{\calO_X} \calQ^j_l,$$ where $(\calK^{(2)})^i_k \otimes_{\calO_X} \calQ^j_l$ is in cohomological degree $i+j$ (and $k$ is even if the term is non-zero). Again the differential is the sum of four terms $d_1:=d_{\calK^{(2)}} \otimes \Id_{\calQ}$, $d_2=\Id_{\calK^{(2)}} \otimes d_{\calQ}$, $d_3$ the Koszul differential induced by the action of $\calV$ and $\calV^{\vee}$, and $d_4$ the Koszul differential induced by the action of $\calW$ and $\calW^{\vee}$. The effect of these differentials on the indices $i,j,k,l$ can be described as follows: $$d_1 : i \to i+1, \quad d_2 : j \mapsto j+1, \quad d_3: \left\{ \begin{array}{ccc} i & \mapsto & i+2 \\ j & \mapsto & j-1 \\ k & \mapsto & k-2 \\ l & \mapsto & l+2 \end{array} \right. , \quad d_4: \left\{ \begin{array}{ccc} i & \mapsto & i+1 \\ k & \mapsto & k-2 \\ l & \mapsto & l+2 \end{array} \right. .$$ One has $$(d_1 + d_2)^2=0, \qquad (d_3 + d_4)^2=0.$$ Hence, disregarding the internal grading, $\mathscr{A} \circ \mathscr{B}(\calQ)$ is the total complex of the double complex with $(p,q)$-term $$C^{p,q}:=\bigoplus_{\genfrac{}{}{0pt}{}{p=-l-3k/2,}{q=i+j+l+3k/2}} (\calK^{(2)})^i_k \otimes_{\calO_X} \calQ^j_l,$$ and with differentials $d'=d_3 + d_4$, $d''=d_1 + d_2$. We know that $\calQ^j_l=0$ if $j+l \ll 0$. Moreover, one checks easily that $(\calK^{(2)})^i_k=0$ if $i+3k/2 \ll 0$. Hence $C^{p,q}=0$ if  $q \ll 0$. It follows, by Proposition \ref{prop:spectralsequence}, that there is a converging spectral sequence $$E_1^{p,q}=H^q(C^{p,*},d'') \Rightarrow H^{p+q}(\mathscr{A} \circ \mathscr{B}(\calQ)).$$ Similarly, disregarding the internal grading, $\calQ$ is the total complex of a double complex $C'$, and $\psi$ is induced by a morphism of double complexes $C' \to C$. Hence it is enough to prove that the morphism induced by $\psi$ from $\calQ$ to $\calK^{(2)} \otimes_{\calO_X} \calQ$, endowed with the differential $d_1 + d
_2$, is a quasi-isomorphism.

Once more, this follows from a spectral sequence argument, using the property that $H(\calK^{(2)})=\calO_X$ (see Lemma \ref{lem:Koszul2}).\end{proof}

\subsection{Regrading} \label{ss:regrading}

In this subsection we introduce a ``regrading'' functor. This functor will play a technical role in \ref{ss:finiteness}, and a more crucial role later in the geometric interpretation of the equivalence.

Consider the functor $$\xi : \calC_{\gr}(\calS) \to \calC_{\gr}(\calR)$$ which sends the $\calS$-dg-module $\calM$ to the $\calR$-dg-module with $(i,j)$-component $\xi(\calM)^i_j:=\calM^{i-j}_j$, the differential and the $\calR$-action on $\xi(\calM)$ being induced by the differential and the $\calS$-action on $\calM$. This functor is clearly an equivalence of categories, and it induces equivalences, still denoted $\xi$, $$\calC^{\searrow}_{\gr}(\calS) \xrightarrow{\sim} \calC^{\swarrow}_{\gr}(\calR), \qquad \calD^{\searrow}_{\gr}(\calS) \xrightarrow{\sim} \calD^{\swarrow}_{\gr}(\calR).$$

\subsection{Categories with finiteness conditions} \label{ss:finiteness}

In the rest of this section we prove that the equivalences $\overline{\mathscr{A}}$ and $\overline{\mathscr{B}}$ restrict to equivalences between subcategories of dg-modules whose cohomology is locally finitely generated. This will eventually allow us to get rid of the technical conditions ``$\nwarrow$'' and ``$\searrow$''.

Let us introduce some more notation. If $\calA$ is a quasi-coherent $\Gm$-equivariant dg-algebra, and if $*=\nwarrow, \swarrow, \searrow, \nearrow$, we denote by $\calC^{*,\fg}_{\gr}(\calA)$, respectively $\calD^{*,\fg}_{\gr}(\calA)$, the full subcategory of $\calC^*_{\gr}(\calA)$, respectively $\calD^*_{\gr}(\calA)$, whose objects are the dg-modules $\calM$ such that $H(\calM)$ is a locally finitely generated $H(\calA)$-module.

We also denote by $\calC \calF \calG_{\gr}(\calA)$ the full subcategory of $\calC_{\gr}(\calA)$ whose objects are the locally finitely generated $\Gm$-equivariant $\calA$-dg-modules, and by $\calD \calF \calG_{\gr} (\calA)$ the localization of the homotopy category of $\calC \calF \calG_{\gr}(\calA)$ with respect to quasi-isomorphisms. Finally we denote by $\calD^{\fg}_{\gr}(\calA)$ the full subcategory of $\calD_{\gr}(\calA)$ whose objects are the $\Gm$-equivariant dg-modules $\calM$ such that $H(\calM)$ is locally finitely generated over $H(\calA)$. 

We are going to prove that, in the cases we are interested in, several of these categories coincide. Observe in particular that there are inclusions $$\calC \calF \calG_{\gr}(\calR) \hookrightarrow \calC^{\swarrow,\fg}_{\gr}(\calR), \quad \calC \calF \calG_{\gr}(\calS) \hookrightarrow \calC^{\searrow,\fg}_{\gr}(\calS), \quad \calC \calF \calG_{\gr}(\calT) \hookrightarrow \calC^{\nwarrow,\fg}_{\gr}(\calT),$$ which induce functors between the corresponding derived categories.

\begin{lem} \label{lem:equivalencesfg}

$\rmi$ The induced functors $$\calD \calF \calG_{\gr}(\calR) \to \calD^{\swarrow,\fg}_{\gr}(\calR), \ \ \calD \calF \calG_{\gr}(\calS) \to \calD^{\searrow,\fg}_{\gr}(\calS), \ \ \calD \calF \calG_{\gr}(\calT) \to \calD^{\nwarrow,\fg}_{\gr}(\calT)$$ are equivalences of categories.

$\rmii$ Similarly, the natural functors $$\calD \calF \calG_{\gr}(\calR) \to \calD^{\fg}_{\gr}(\calR), \quad \calD \calF \calG_{\gr}(\calS) \to \calD^{\fg}_{\gr}(\calS), \quad \calD \calF \calG_{\gr}(\calT) \to \calD^{\fg}_{\gr}(\calT)$$ are equivalences of categories.

\end{lem}

\begin{proof} Our proof of this lemma is very similar to that of \cite[VI.2.11]{BDmod} (see also \cite[3.3.4]{RKos}). We give the details of the proof of $\rmii$. Statement $\rmi$ can be treated similarly.

Using the ``regrading trick'' of \ref{ss:regrading}, the cases of $\calS$ and $\calR$ are equivalent. Similarly, using the change of the internal grading to the opposite one, we see that the cases of $\calR$ and $\calT$ are equivalent. Hence it is sufficient to consider the $\Gm$-equivariant dg-algebra $\calT$. 

Remark that the algebra $\calT$, as well as its cohomology $H(\calT)$, is finitely generated as a ${\rm S}(\calW)$-module. Hence a $\calT$-dg-module $\calN$ is locally finitely generated, respectively has locally finitely generated cohomology, iff $\calN$, respectively $H(\calN)$, is locally finitely generated over ${\rm S}(\calW)$.

\begin{lem} \label{lem:inductivelimit}

Let $\calN$ be an object of $\calC_{\gr}(\calT)$, with locally finitely generated cohomology, whose cohomological grading is bounded. Then $\calN$ is the inductive limit of quasi-coherent sub-$\calT$-dg-modules which are locally finitely generated, and which are quasi-isomorphic to $\calN$ under inclusion.

\end{lem}

\begin{proof}[Proof of Lemma \ref{lem:inductivelimit}] The internal grading has no importance in this statement, hence we will forget about it in the proof. The dg-module $\calN$ is clearly an inductive limit of locally finitely generated quasi-coherent sub-$\calT$-dg-modules. Hence it is sufficient to show that given a locally finitely generated quasi-coherent sub-dg-module $\calF$ of $\calN$, there exists a locally finitely generated quasi-coherent sub-dg-module $\calG$ of $\calN$, containing $\calF$ and
quasi-isomorphic to $\calN$ under the inclusion map.

This is proved by a simple (descending) induction. Let $j \in \mathbb{Z}$. Assume that we have found a subcomplex $\calG_{(j)}$ of $\bigoplus_{i \geq j} \calN^i$, quasi-coherent over $\calO_X$, locally finitely generated over ${\rm S}(\calW)$, containing $\bigoplus_{i \geq j} \calF^i$, stable under $\calT$ (\emph{i.e.} if $g \in \calG_{(j)}^i$ and $t \in \calT^k$, and if $i+k \geq j$, then $t \cdot g \in \calG_{(j)}^{i+k}$), such that $\calG_{(j)} \hookrightarrow \calN$ is a quasi-isomorphism in cohomological degrees greater than $j$ and that $\calG^j_{(j)} \cap \ker(d_{\calN}^j) \to H^j(\calN)$ is surjective. Then we choose a locally finitely generated sub-${\rm S}(\calW)$-module $\calH^{j-1}$ of $\calN^{j-1}$ containing $\calF^{j-1}$, quasi-coherent over $\calO_X$, whose image under $d_{\calN}^{j-1}$ is $\calG^j_{(j)} \cap \Ima(d_{\calN}^{j-1})$. Without altering these conditions, we can add a sub-module of cocycles so that the new sub-module
$\calH^{j-1}$ contains representatives of all the elements of
$H^{j-1}(\calN)$. We can also assume that $\calN^{j-1}$ contains all the
sections of the form $t \cdot g$ for $t \in \calT^i$ and $g \in
\calG_{(j)}^k$ with $i+k=j-1$. Then we define $\calG_{(j-1)}$ by
$$\calG_{(j-1)}^k=\left\{ \begin{array}{cl} \calG_{(j)}^k & \text{if } k
    \geq j, \\ \calH^{j-1} & \text{if } k=j-1. \end{array} \right.$$
For $j$ small enough, $\calG_{(j)}$ is the sought-for sub-dg-module. \end{proof}

Let us denote by $$\iota : \calD \calF \calG_{\gr}(\calT) \to \calD_{\gr}^{\fg}(\calT)$$ the functor under consideration. Let $\calN$ be an object of $\calD^{\fg}_{\gr}(\calT)$. Then the cohomology $H(\calN)$ is bounded for the cohomological grading (because it is locally finitely generated over $H(\calT)$, which is bounded). Hence, using truncation functors (see \ref{ss:notations}), $\calN$ is isomorphic to a $\calT$-dg-module whose cohomological grading is bounded. Using Lemma \ref{lem:inductivelimit}, it follows that $\calN$ is in the essential image of $\iota$. Hence $\iota$ is essentially surjective.

Now, let us prove that it is full. Let $\calN_1$ and $\calN_2$ be objects of $\calC \calF \calG_{\gr}(\calT)$. In particular, $\calN_1$ and $\calN_2$ have bounded cohomological grading. A morphism $\phi: \iota(\calN_1) \to \iota(\calN_2)$ in $\calD^{\fg}_{\gr}(\calT)$ is represented by a diagram $$\iota(\calN_1) \xrightarrow{\alpha} \calF \xleftarrow{\beta} \iota(\calN_2)$$ where $\beta$ is a quasi-isomorphism. Using truncation functors, one can assume that $\calF$ has bounded cohomological grading. By Lemma \ref{lem:inductivelimit}, there exists a locally finitely generated sub-$\calT$-dg-module $\calF '$ of $\calF$, containing $\alpha(\calN_1)$ and $\beta(\calN_2)$, and quasi-isomorphic to $\calF$ under the inclusion map. Then $\phi$ is also represented by $$\iota(\calN_1) \xrightarrow{\alpha} \calF ' \xleftarrow{\beta} \iota(\calN_2),$$ which is the image of a morphism in $\calD \calF \calG_{\gr}(\calT)$. Hence $\iota$ is full.

Finally we prove that $\iota$ is faithful. If a morphism $f:\calN_1 \to
\calN_2$ in $\calC \calF \calG_{\gr}(\calT)$ is such that $\iota(f)=0$, then
there exist an object $\calF$ of $\calD^{\fg}_{\gr}(\calT)$, which can again
be assumed to be bounded, and a quasi-isomorphism of $\calT$-dg-modules $g: \calN_2 \to \calF$ such that $g \circ f$ is homotopic to zero. This homotopy is given by a morphism $h: \calN_1 \to \calF[-1]$. Again by Lemma \ref{lem:inductivelimit}, there exists a locally finitely generated sub-$\calT$-dg-module $\calF '$ of $\calF$ containing $g(\calN_2)$ and $h(\calN_1)[1]$, and quasi-isomorphic to $\calF$ under inclusion. Replacing $\calF$ by $\calF'$, this proves that $f=0$ in $\calD \calF \calG_{\gr}(\calT)$. The proof of Lemma \ref{lem:equivalencesfg} is complete. \end{proof}

\subsection{Restriction of the equivalences to locally finitely generated dg-mo\-dules} \label{ss:restrictionfiniteness}

\begin{prop} \label{prop:restrictionequiv}

The functors $\overline{\mathscr{A}}$, $\overline{\mathscr{B}}$ restrict to equivalences of categories $$\calD^{\searrow,\fg}_{\gr}(\calS) \cong \calD^{\nwarrow,\fg}_{\gr}(\calT).$$

\end{prop}

\begin{proof} It is sufficient to prove that the functors $\overline{\mathscr{A}}$, $\overline{\mathscr{B}}$ send dg-modules with locally finitely generated cohomology to dg-modules with locally finitely generated cohomology.

\emph{First step: functor $\overline{\mathscr{B}}$.} First, let us consider $\overline{\mathscr{B}}$. By Lemma \ref{lem:equivalencesfg}, it suffices to prove that if $\calN$ is a locally finitely generated $\calT$-module, then $\overline{\mathscr{B}}(\calN)$ has locally finitely generated cohomology. We begin with the following lemma.

\begin{lem} \label{lem:resolution2}

Let $\calN$ be a locally finitely generated $\Gm$-equivariant $\calT$-dg-module. There exist an object $\calQ$ of $\calC \calF \calG_{\gr}(\calT)$, which is locally free of finite rank over ${\rm S}(\calW) \subset \calT$, and a quasi-isomorphism $\calQ \xrightarrow{\qis} \calN$.

\end{lem}

\begin{proof}[Proof of Lemma \ref{lem:resolution2}] The arguments in this proof are very close to those in the proof of Proposition \ref{prop:resolutions}. There exists a $\Gm$-equivariant sub-$\calO_X$-dg-module $\calG \subset \calN$, which is coherent as an $\calO_X$-module, and which generates $\calN$ as a $\calS$-dg-module. There exists also a $\Gm$-equivariant $\calO_X$-dg-module $\calF$, which is locally free of finite rank as an $\calO_X$-module, and a surjection $\calF \twoheadrightarrow \calG$. We set $$\calQ^{(1)}:=\calT \otimes_{\calO_X} \calF,$$ endowed with its natural structure of $\Gm$-equivariant $\calT$-dg-module. Then we have a surjection of $\calT$-dg-modules $$\calQ^{(1)} \twoheadrightarrow \calN,$$ and $\calQ^{(1)}$ is locally free over ${\rm S}(\calW)$.

Let $n$ be the rank of $\calW$ over $\calO_X$. Taking the kernel of our morphism $\calQ^{(1)} \to \calN$, and repeating the argument, we obtain locally finitely generated $\calT$-dg-modules $\calQ^{(j)}$, $j=1, \cdots, n+d$, which are locally free of finite rank over ${\rm S}(\calW)$, and an exact sequence of $\calT$-dg-modules $$\calQ^{(n+d)} \to \calQ^{(d+n-1)} \to \cdots \to \calQ^{(1)} \to \calN \to 0.$$ All these objects are complexes of coherent $\calS(\calW)$-modules, hence we can consider them as complexes of coherent sheaves on $W^*$, the vector bundle on $X$ with sheaf of sections $\calW^{\vee}$. The scheme $W^*$ is noetherian, integral, separated and regular of dimension $d+n$. Hence $\calQ^{(n+d+1)}:={\rm Ker}(\calQ^{(n+d)} \to \calQ^{(n+d-1)})$ is also locally free over ${\rm S}(\calW)$. Then $$\calQ:={\rm Tot}(0 \to \calQ^{(n+d+1)} \to \cdots \to \calQ^{(1)} \to 0)$$ is a resolution of $\calN$ as in the lemma.\end{proof}

Now let $\calQ \xrightarrow{\qis} \calN$ be a resolution as in Lemma \ref{lem:resolution2}. In particular $\calQ$ is locally free over $\calO_X$, hence $\overline{\mathscr{B}}(\calN)$ is isomorphic to the image of $\mathscr{B}(\calQ)$ in the derived category. Hence it is enough to prove that $\mathscr{B}(\calQ)$ has locally finitely generated cohomology, and even to prove that this cohomology is locally finitely generated over ${\rm S}(\calV^{\vee})$. Let $x \in X$. The object $\mathscr{B}(\calQ)_x$ was described in \ref{ss:derivedfunctors}. We use the same notations as in this subsection. Disregarding the internal grading, $\mathscr{B}(\calQ)_x$ is also the total complex of the double complex with $(p,q)$-term $$C^{p,q}:=\bigoplus_{\genfrac{}{}{0pt}{}{p=j,}{q=i+k+j}} \Lambda^i(\calW_x^{\vee}) \otimes_{\calO_{X,x}} {\rm S}^j(\calV_x^{\vee}) \otimes_{\calO_{X,x}} (\calQ_x^{\vee})^k_l,$$ and with differentials $d'=d_1+d_3$, $d''=d_2+d_4$. By hypothesis, $(\calQ^{\vee})^k_l=0$ for $k \ll 0$, hence $C^{p,q}=0$ for $q \ll 0$. Hence by Proposition \ref{prop:spectralsequence} there is a converging spectral sequence $$E_1^{p,q}=H(C^{p,*},d'') \Rightarrow H^{p+q}(\mathscr{B}(\calQ)_x).$$ It follows that it is sufficient to prove that the cohomology of $\calS \otimes_{\calO_X} \calQ^{\vee}$, endowed with the differential $d_2+d_4$, is locally finitely generated over ${\rm S}(\calV^{\vee})$. This complex is again the total complex of the double complex with $(p,q)$-term $$D^{p,q}:=\bigoplus_{\genfrac{}{}{0pt}{}{p=2j+k,}{q=i}} \Lambda^i(\calW_x^{\vee}) \otimes_{\calO_{X,x}} {\rm S}^j(\calV_x^{\vee}) \otimes_{\calO_{X,x}} (\calQ_x^{\vee})^k_l,$$ and with differentials $d'=d_2$, $d''=d_3$. The spectral sequence of this double complex again converges, hence we can forget about $d_2$. Then $\calS \otimes_{\calO_X} \calQ^{\vee}$, endowed with the differential $d_3$, is locally the tensor product of ${\rm S}(\calV^{\vee})$ with a finite number of Koszul complexes ${\rm Koszul}_2(\calW_x^{\vee})$ of \eqref{eq:Koszulcomplex2}. The result follows. \medskip

\emph{Second step: functor $\overline{\mathscr{A}}$.} The proof for the functor $\overline{\mathscr{A}}$ is entirely similar. In this case, with the notation of \ref{ss:derivedfunctors}, we can consider the double complex with $(p,q)$-term $$C^{p,q}:=\bigoplus_{\genfrac{}{}{0pt}{}{p=k-2i-j,}{q=i+j}} \Lambda^i(\calW_x^{\vee}) \otimes_{\calO_{X,x}} {\rm S}^j(\calV_x^{\vee}) \otimes_{\calO_{X,x}} (\calP_x^{\vee})^k_l,$$ and differentials $d'=d_1+d_2$, $d''=d_3+d_4$. Here $C^{p,q}=0$ for $q<0$, hence the corresponding spectral sequence converges, and we can forget about $d_1$ and $d_2$. Then we can consider the double complex \[ D^{p,q}:=\bigoplus_{\genfrac{}{}{0pt}{}{p=k-2i,}{q=i}} \Lambda^i(\calW_x^{\vee}) \otimes_{\calO_{X,x}} {\rm S}^j(\calV_x^{\vee}) \otimes_{\calO_{X,x}} (\calP_x^{\vee})^k_l,\] with differentials $d'=d_4$ and $d''=d_3$. And we finish the proof as above. \end{proof}

Finally, combining Proposition \ref{prop:restrictionequiv}, Lemma \ref{lem:equivalencesfg} and the ``regrading trick'' of \ref{ss:regrading} we obtain the following theorem, which is the main result of this section.

\begin{thm} \label{thm:algduality}

There exists a contravariant equivalence of triangulated categories $$\kappa: \calD^{\fg}_{\gr}(\calT) \xrightarrow{\sim} \calD^{\fg}_{\gr}(\calR)$$ satisfying $\kappa(\calM[n] \langle m \rangle)=\kappa(\calM)[-n+m] \langle m \rangle$.

\end{thm}

\section{Linear Koszul duality}

In this section we give a geometric interpretation of Theorem \ref{thm:algduality}.

\subsection{Intersections of vector bundles}

Let us consider as above a noetherian, integral, separated, regular scheme $X$, and a vector bundle $E$ over $X$. Let $F_1, F_2 \subset E$ be sub-vector bundles. Let $E^*$ be the vector bundle dual to $E$, and let $F_1^{\bot}, F_2^{\bot} \subset E^*$ be the orthogonal to $F_1$, respectively $F_2$. We will be interested in the dg-schemes $$F_1 \, \rcap_E \, F_2 \qquad \text{ and } \qquad F_1^{\bot} \, \rcap_{E^*} \, F_2^{\bot}.$$

Let $\calE, \calF_1, \calF_2$ be the sheaves of local sections of $E, F_1, F_2$. Then the sheaves of local sections of $E^*, F_1^{\bot}, F_2^{\bot}$ are, respectively, $\calE^{\vee}$, $\calF_1^{\bot}$ and $\calF_2^{\bot}$ (here we consider the orthogonals inside $\calE^{\vee}$). Let us denote by $\calX$ the $\calO_X$-dg-module $$\calX:= \bigl( 0 \to \calF_1^{\bot} \to \calF_2^{\vee} \to 0 \bigr),$$ where $\calF_1^{\bot}$ is in degree $-1$, $\calF_2^{\vee}$ is in degree $0$, and the non-trivial differential is the composition of the natural morphisms $\calF_1^{\bot} \hookrightarrow \calE^{\vee} \twoheadrightarrow \calF_2^{\vee}$, and by $\calY$ the $\calO_X$-dg-module \[ \calY:= \bigl( 0 \to \calF_2 \to \calE / \calF_1 \to 0 \bigr),\] where $\calF_2$ is in degree $-1$, $\calE / \calF_1$ is in degree $0$, and the non-trivial differential is the opposite of the composition of the natural morphisms $\calF_2 \hookrightarrow \calE \twoheadrightarrow  \calE / \calF_1$.

\begin{lem}

There exist equivalences of categories $$\calD(F_1 \, \rcap_E \, F_2) \cong \calD(X, \Sym(\calX)), \qquad \calD(F_1^{\bot} \, \rcap_{E^*} \, F_2^{\bot}) \cong \calD(X, \Sym(\calY)).$$

\end{lem}

\begin{proof} We need only prove the first equivalence (the second one is similar: replace $E$ by $E^*$, $F_1$ by $F_2^{\bot}$, $F_2$ by $F_1^{\bot}$). Let $\calA$ be any graded-commutative, non-positively graded, quasi-coherent dg-algebra on $E$, quasi-isomorphic to $\calO_{F_1} \, \lotimes_{\calO_E} \, \calO_{F_2}$ (see \ref{ss:derivedintersection}). Let $\pi: E \to X$ be the natural projection. Then it is well-known (see \emph{e.g.} \cite[1.4.3]{EGA2}) that the functor $\pi_*$ induces equivalences of categories $$\calC(E,\calA) \cong \calC(X,\pi_* \calA), \quad \calD(E,\calA) \cong \calD(X,\pi_* \calA).$$ Moreover, the data of $\calA$ is equivalent to the data of the $\pi_* \calO_E$-dg-algebra $\pi_* \calA$, which is quasi-isomorphic to $\pi_* \calO_{F_1} \, \lotimes_{\pi_* \calO_E} \, \pi_* \calO_{F_2}$.

Now there are natural isomorphisms $\pi_* \calO_E \cong {\rm S}_{\calO_X}(\calE^{\vee})$, $\pi_* \calO_{F_i} \cong {\rm S}_{\calO_X}(\calF_i^{\vee})$ ($i=1,2$). Consider the Koszul resolution $$\Sym \bigl( 0 \to \calF_1^{\bot} \to \calE^{\vee} \to 0 \bigr) \ \xrightarrow{\qis} \ {\rm S}(\calF_1^{\vee}) \cong {\rm S}(\calE^{\vee})/(\calF_1^{\bot} \cdot {\rm S}(\calE^{\vee})),$$ where $\calF_1^{\bot}$ is in degree $-1$, $\calE^{\vee}$ is in degree $0$, and the differential is the natural inclusion. This is a flat dg-algebra resolution of ${\rm S}(\calF_1^{\vee})$ over ${\rm S}(\calE^{\vee})$. If we tensor this resolution with ${\rm S}(\calF_2^{\vee})$ (over ${\rm S}(\calE^{\vee})$) we obtain that the dg-algebra $\Sym(\calX)$ is quasi-isomorphic to $\pi_* \calO_{F_1} \, \lotimes_{\pi_* \calO_E} \, \pi_* \calO_{F_2}$. Hence we can take $\pi_* \calA = \Sym(\calX)$. This finishes the proof of the lemma.\end{proof}

\subsection{Linear Koszul duality}

One can also consider $\calX$ as a $\Gm$-equivariant $\calO_X$-dg-module, where $\calF_1^{\bot}$ and $\calF_2^{\vee}$ are in internal degree $2$. Then, similarly, $\calY$ is $\Gm$-equivariant (with generators in internal degree $-2$). Geometrically, this corresponds to considering a $\Gm$-action on $E$, where $t \in \bk^{\times}$ acts by multiplication by $t^{-2}$ along the fibers. We will use the notations \begin{eqnarray*} \calD_{\Gm}^{{\rm c}}(F_1 \, \rcap_E \, F_2) & := & \calD_{\gr}^{\fg}(X, \Sym(\calX)), \\ \calD_{\Gm}^{{\rm c}}(F_1^{\bot} \, \rcap_{E^*} \, F_2^{\bot}) & := & \calD_{\gr}^{\fg}(X, \Sym(\calY)). \end{eqnarray*} Then Theorem \ref{thm:algduality} gives, in this situation:

\begin{thm}

There exists a contravariant equivalence of triangulated categories, called \emph{linear Kos\-zul duality}, $$\kappa: \calD_{\Gm}^{{\rm c}}(F_1 \, \rcap_E \, F_2) \ \xrightarrow{\sim} \ \calD_{\Gm}^{{\rm c}}(F_1^{\bot} \, \rcap_{E^*} \, F_2^{\bot})$$ satisfying $\kappa(\calM[n] \langle m \rangle)=\kappa(\calM)[-n+m] \langle m \rangle$.

\end{thm}

\subsection{
Equivariant version of the duality
}

Finally, let us consider an algebraic group $G$ acting on $X$ (algebraically). Assume that $E$ is a $G$-equivariant vector bundle, and that $F_1$ and $F_2$ are $G$-equivariant subbundles. Then, with the same notations as above, $\calX$ is a complex of $G$-equivariant coherent sheaves on $X$. Let us denote by $$\calD_{G \times \Gm}^{{\rm c}}(F_1 \, \rcap_E \, F_2)$$ the derived category of $G \times \Gm$-equivariant quasi-coherent $\Sym(\calX)$-dg-modules on $X$ (\emph{i.e.} $\Gm$-equivariant dg-modules as above, endowed with a structure of $G$-equivariant quasi-coherent $\calO_X$-module compatible with all other structures) with locally finitely generated cohomology, and similarly for $\calD_{G \times \Gm}^{{\rm c}}(F_1^{\bot} \, \rcap_{E^*} \, F_2^{\bot})$. Then our constructions easily extend to give the following result.

\begin{thm}

There exists a contravariant equivalence of categories $$\kappa: \calD_{G \times \Gm}^{{\rm c}}(F_1 \, \rcap_E \, F_2) \ \xrightarrow{\sim} \ \calD_{G \times \Gm}^{{\rm c}}(F_1^{\bot} \, \rcap_{E^*} \, F_2^{\bot})$$ satisfying $\kappa(\calM[n] \langle m \rangle)=\kappa(\calM)[-n+m] \langle m \rangle$.

\end{thm}

\end{document}